\documentclass[10pt,a4paper,twoside,leqno]{article}
\usepackage{amsfonts,amssymb,amscd}

\usepackage{amsmath,graphicx}
\usepackage[T1]{fontenc}
\usepackage[all]{xy}
\usepackage[cp1250]{inputenc}

\font\Bigtit=cmr10 scaled \magstep 4
\font\ebf=cmbx8

\newtheorem{thm}{Theorem}
\newtheorem{lem}{Lemma}
\newtheorem{prop}{Proposition}

\renewcommand{\H}{\mathcal H}

\newcommand{\CC}{\mathcal C}

\newcommand{\x}{\widetilde x}
\newcommand{\w}{\widetilde w}
\newcommand{\y}{\widetilde y}
\newcommand{\X}{\widetilde X}

\newcommand{\U}{\widetilde U}
\newcommand{\V}{\widetilde V}
\newcommand{\f}{\widetilde f}

\newcommand{\F}{\mathcal F}

\newcommand{\Tomega}{\widetilde \omega}

\newcommand{\TA}{\widetilde{\mathcal A}}

\newcommand{\al}{\alpha}
\newcommand{\be}{\beta}

\newcommand{\FF}{\mathcal F}
\newcommand{\A}{\mathcal A}

\newcommand{\R}{\mathbb R}
\newcommand{\Z}{\mathbb Z}
\newcommand{\N}{\mathbb N}

\newenvironment{Proof of}[1]{\textbf{Proof of the #1.}}{$\qquad \blacksquare$\par}
\newenvironment{Item}[1]{\par #1 }{\par}

\begin{document}

\thispagestyle{empty}

\begin{flushright}
PL ISSN 0459-6854
\end{flushright}
\vspace{0.5cm}
\centerline{\Bigtit B U L L E T I N}
\vspace{0.5cm}
\centerline{DE \ \  LA \ \  SOCI\'ET\'E \ \  DES \ \  SCIENCES \ \ ET \ \ DES \
\ \ LETTRES \ \ DE \ \ \L \'OD\'Z}
\vspace{0.3cm}
\noindent 2005\hfill Vol. LV
\vspace{0.3cm}
\hrule
\vspace{5pt}
\noindent Recherches sur les d\'eformations \hfill Vol. XLVIII
\vspace{5pt}
\hrule
\vspace{0.3cm}
\noindent pp.~83--109

\vspace{0.7cm}

\noindent {\it Bartosz K.  Kwa\'sniewski}

\vspace{0.5cm}

\noindent {\bf  INVERSE LIMIT SYSTEMS ASSOCIATED\\ WITH
$\mathcal{F}_{2^n}$ ZERO SCHWARZIAN UNIMODAL  MAPS}

\vspace{0.5cm}

\noindent {\ebf Summary}

{\small We present an illustrative example of  an  inverse limit space and
a shift map associated with an ${\F}_{2^n}$ unimodal mapping
consisting of two hyperbolae. Topologically, in case $n=0$  the
limit space is an interval, in case $n=1,2$, it is a $\sin
\frac{1}{x}$-continuum, and in case $n=3$ it is a certain
continuum   endowed with a specific geometrical beauty. The
dynamics of the shift map is also described.

The article is motivated by the recent discovery \cite{maxid},
\cite{covalg} that inverse limit structures play  fundamental role in
crossed-products of $C^*$-algebras, and hence  this paper is possibly of
interest not only to topologists or dynamicists but also to operator
algebraists.}

\vspace{0.5cm}

%section 1
\subsection{Introduction}

Inverse limit\footnote{another name for inverse limit is   projective limit} spaces received much attention by
topologists and   dynamicists for its application in continuum theory
\cite{Nadler}, \cite{Watkins}, \cite{Brucks Bruin}, hyperbolic attractors
\cite{Williams1}, \cite{Williams3}, \cite{Williams2}, \cite{Yi}, and
recently  substitution tilings \cite{Anderson Putnam}. Usually, the advantage of the
inverse limits  is that they reduce investigation of 'large'
topologically complicated objects to much 'smaller' and consequently nicer ones. Thus in this approach  the inverse limit spaces are initial objects. However, the quite
opposite way of thinking proved also to be useful for instance in the
theory of wavelets and martingales \cite{Dutkay Jorgensen}, and in  crossed
products of $C^*$-algebras (see section \ref{C^*} for a brief survey on
this  context) where the inverse limit system (which is always reversible)
helps to understand the initial  irreversible system. 

Apart from the constant interest  on the subject and many
beautiful results, inverse limit spaces still seem to be abstract and
not completely accessible.   There are  not many examples   of the explicitly calculated  inverse limit spaces and
their shift maps in the
literature.   The model examples are Knaster continua
and $n$-adic solenoids, see \cite{Nadler}, \cite{Watkins},
\cite{introdynsys}.

The aim of the present paper is to provide   an example
which would accustom even non-specialists to inverse limits, and  which
would give the opportunity to see  (literally) how the topology of inverse
limit space and dynamics of the shift map depend on the dynamics of the
initial map.

The paper is organized as follows. In Section\,2 we recall
the definition of an inverse limit space $\X$ and its shift map
$\f$ induced by a single bonding map $f:X\to X$, we introduce  the
class of bonding maps which we divide (depending on the complexity
of dynamics) into the following subclasses 1), 2) 3a) and 3b);
then we briefly present  the main results. Section 3
presents the methods  we  used to
achieve the goal while Sections\,4 and 5 contain
detailed proves and results for cases 1), 2) and 3a), 3b)
respectively. Lastly, in Section 6 we introduce the interpretation
of the inverse limit system $(\X,\f)$ in terms of $C^*$-algebras and
explain their applications to crossed-products.

\subsection{The basic object of the paper and the main result}\label{2}

Let  $X$ be a compact topological space and $f:X\rightarrow X$  a continuous
surjective but noninjective mapping.
Then $(X,f)$ is a typical  irreversible topological dynamical system and there
is a natural way of constructing an
 extended reversible dynamical system that is a pair $(\X,\f)$ where $\X$ is
 compact, $\f$ is a homeomorphism,
  and there exists a continuous surjective mapping $\Phi:\X\rightarrow X$
  such that the following diagram commutes
\begin{equation}\label{phidiagram}
 \begin{array}{c c c} \X & \stackrel{\f}{\longrightarrow} & \X  \\
\Phi\downarrow &      & \downarrow \Phi\\
X & \stackrel{f}{\longrightarrow} & X,
\end{array}\end{equation}

\noindent hence $(X,f)$ is a factor of $(\X,\f)$. It is enough to  take
  $\X$ to be the inverse limit of the inverse sequence $X\stackrel{f}{\longleftarrow}X
\stackrel{f}{\longleftarrow} ...$, and $\f$ to be the mapping induced by $f$. Namely,
\begin{equation}\label{inverse space}
\X=\{(x_0,x_1,...)\in X^\N: f(x_{n+1})=x_{n},\,\, n\in \N\}
\end{equation}

\noindent (we adhere to the convention that natural numbers starts from zero) is
 equipped with the topology inherited from the product space $X^\N$, and $\f$ acts as follows
\begin{equation}\label{inverse mapping}
\f(x_0,x_1,x_2,...)=(f(x_0),x_0,x_1,...), \qquad \textrm{ for }(x_0,x_1,...)\in \X.
\end{equation}

\noindent The factor map $\Phi:\X\rightarrow X$ is then given by
$$
\Phi(x_0,x_1,x_2,...)=x_0.
$$

\noindent We shall call the pair $(\X,\f)$, constructed in the above manner,
a \textsl{reversible extension } of  $(X,f)$, cf. \cite{covalg}.

It is clear that in order to obtain a relatively nice representation of an inverse
limit space $\X$
 (a rather complicated object)  an initial  system $(X,f)$ should be relatively
 simple.
Our choice of $(X,f)$ was dictated by the following requirements:
\begin{itemize}
\item[1)] The dynamics should be low-dimensional
\item[2)] Mapping $f$ should be piece-wise monotone
\item[3)] The dynamical system  should not be too chaotic
\end{itemize}
According to 1) and 2) we choose $X$ to be an interval $[0,1]$,
and $f$ to be a unimodal mapping, that  is a piece-wise monotone
mapping possessing  one extremum.

 According to 3) we impose on $f$ the condition to have finitely many periodic points.
Then by Sharkovsky's theorem (see, for instance, \cite{introdynsys}) there exist
$N\in\N$ such that $Fix(f^{2^N})=Fix(f^{2^{N+1}})$
where $Fix(g)$ denotes the set of fixed points of mapping $g$.
Let us recall that $f$ is said to be $\F_{2^N}$ if $N$ is the smallest number
with the above described property.

As even non-experts nowadays know, see \cite{Stewart}, already
such a simple mapping as quadratic one may  ``cause chaos''.
Hence  the third requirement makes our choice not an easy task. However, the paper
of A.\,Bondarenko and  S.\,Popovych \cite{BonPop} help us out. We borrow an
appropriate  mapping from them.

A.\,Bondarenko and  S.\,Popovych reduced  investigation of a
continuous unimodal map $f:[0,1]\rightarrow [0,1]$  of the form:
$$
f(x)=\left \{ \begin{array}{ll} f_0(x)=\displaystyle{\frac{\al_1 x+\be_1}{\gamma_1 x+
\delta_1}} , & x\in [0,\rho],\\
f_1(x)=\displaystyle{\frac{\al_2 x+\be_2}{\gamma_2 x+ \delta_2}} , & x\in
(\rho,1],\end{array} \right.
$$

\noindent to the investigation of maps of two types, and here we restrict ourselves to
mappings only of type 1, see   \cite{BonPop}.
Namely we assume that $f$ has the following form (see Fig.\,1):
\begin{equation}\label{f(x)}
f(x)=\left \{ \begin{array}{ll} f_0(x)=\displaystyle{\frac{\al x- \al\rho}{\gamma x+
\delta}} , & x\in [0,\rho],\\
f_1(x)=\displaystyle{\frac{x-\rho}{1-\rho}} , & x\in (\rho,1],\end{array} \right.
\end{equation}where parameters $\rho$, $\delta$, $\gamma$, $\al$ range through all
real numbers satisfying the relations
$$
0 < \rho < 1,\qquad 0< \delta,\qquad -\,\frac{\delta}{\rho}< \gamma,\qquad -\,
\frac{\delta}{\rho}< \al<0.
$$

\noindent We denote by $\omega_0$ a fixed point in $[0,\rho)$ and define
$\rho_1$ as $f_0(0)$ (see Fig.\,1).

The  goal of the present paper is a  description of a reversible
extension $(\X,\f)$ of a
dynamical system $(X,f)$  where $X=[0,1]$ and $f$ is an $\F_{2^n}$ mapping
given by (\ref{f(x)}).
We shall obtain it with help of the  following  lemma  borrowed from \cite{BonPop},
see \cite[Lemma\,1]{BonPop}.
\begin{lem}\label{lemma1}
Let $(X,f)$ be $\F_{2^n}$ dynamical system. Then $n\leq 2$ and only following cases
are possible:
\end{lem}
\begin{enumerate}
\item {\it If  $\rho > \rho_1$,  then  $n=0$, $\omega_0$ is attractive and there are no
cycles of larger periods.}
\item {\it If  $\rho = \rho_1$,  then  $n=1$, the left hyperbola is symmetric  with respect
to diagonal.
Therefore each point of $[0,\rho]\setminus\{\omega_0\}$ has period two. }
\item {\it If  $\rho < \rho_1$,  then   $\omega_0$ is repellent and either

{ a)} $n=1$, there exists an attractive cycle of period two and there are no cycles of
larger periods, or

{ b)} $n=2$, there exist two intervals of periodic points such that the middle points
 of those
 intervals have period two and other points have period four
 and  no cycles of larger periods exist. }
\end{enumerate}

\begin{center}
\setlength{\unitlength}{1mm}
\begin{picture}(50,57)(0,-6)
\thinlines
\put(0,0){\framebox(50,50){}}
\qbezier[200](0,0)(25,25)(50,50)
\thicklines
\qbezier(0,27)(20,20)(30,0)
\qbezier(30,0)(40,25)(50,50)
\thinlines
\put(-1,-4){$0$}
\put(29,-3){$\rho$}
\put(49,-4){$1$}
\put(-4.5,26){$\rho_1$}

\qbezier[50](17,0)(17,8.5)(17,17)\put(15.5,-3){$\omega_0$}

\put(-3,49){$1$}
%\put(70,23){\textbf{Figure\,1.}}
\end{picture}

\noindent {\small Fig.\,1. }
\end{center}

The description of $(\X,\f)$ depends on $(X,f)$. If $(X,f)$ falls
into case  1, 2, 3a), or 3b) from the above lemma, then $(\X,\f)$
is described in Theorem \ref{atrractivepointtheorem},
\ref{rho_1=rhotheorem}, \ref{rho_1>rhotheorem}, or
\ref{rho_1>rhotheoremb)} respectively. While viewing  these
results one may find it natural that the induced mapping $\f$
depends  on the number and character of periodic points of $f$,
that means the  homeomorphism $\f$  is basically  different  in cases
1, 2, 3a), 3b). However, one may  also find  it  interesting that
the inverse limit space $\X$ depends only on  the number $n$,
and if $(X,f)$ is $\F_{2^n}$, then $\X$ can be embedded into
$\mathbb{R}^{n+1}$, so $n$ occurs here as  'a dimension' of $\X$. The reader
familiar with the subject will recognize  the typical arc+ray subcontinua \cite{Brucks Bruin}.
Namely, Theorems 2,
3, 4,
5 imply
\begin{thm}
Let $(X,f)$ be as in the preceding lemma. Then depending on $n$,
 the inverse limit space $\X$ can be considered as the following subspace of
 $\mathbb{R}^{n+1}$:

\noindent {\rm 1)} For $n=0$,  $\X$ is an arc, $\X=[0,1]$.

\noindent {\rm 2)} For $n=1$,  $\X$ is the one sided arc+ray continuum, see Figs\,3, 7,
i.e. it is the closure of the graph of an oscillating
function:
$$\X=\underbrace{\{0\}\times[-1,1]}_{arc}\cup \underbrace{\big\{(x,\sin\frac{1}{x}):
x\in (0,1]\big\}}_{ray}.$$

\noindent {\rm 3)} For $n=2$,  $\X$ consists of two arcs and two rays, it  is a two-sided
arc+ray continuum (see Fig.\,8)
with an additional ray winded on it, see Fig.\,9:
$$
\X=\underbrace{\Big\{(x,g(x),0):
x \in (0,1)\Big\}}_{ray} \cup \underbrace{\{0,1\}\times
[-1,1]\times\{0\}}_{two \,\, arcs}\cup \underbrace{\Big\{(\sin\frac{1}{x},
g(\sin\frac{1}{x}),x):x \in (0,1]\Big\}}_{additional\,\, ray}
,$$

\noindent  where
$
g(x)=\left \{ \begin{array}{ll} \sin\frac{\pi}{x} , & x\in (0,\frac{1}{2}],\\
\sin\frac{\pi}{1-x} , & x\in (\frac{1}{2},1)\end{array} \right.
$.

The dynamics of the shift map $\f$ is explained in the Theorems 2,
3, 4 and
5.
\end{thm}

\subsection{The tools and  the  method}\label{3}

We shall build homeomorphic images of the inverse limit space $\X$ using only two simple
lemmas \ref{lem2} and \ref{lem3}   Although they are true in greater generality in order
to avoid unncessary confusion we  assume that
from now on  $X$ will denote the unit interval $[0,1]$,
 $f$ will be an $\F_{2^n}$ mapping given by (\ref{f(x)}), and $(\X,\f)$ will be the
 reversible extension of $(X,f)$.

\noindent At first let us observe that  we can  encode  a point $\x=(x_0,x_1,...)$
from $\X$
by the point $x_0\in X$ and a sequence $T\in  \{0,1\}^{\N}$ such that  $n$-th index
$T(n)$ of $T$ equals to $0$
 iff $x_{n+1}=f^{-1}_0(x_{n})$.
 This code, that is a pair $(x_0,T)$, determines  $\x$ uniquely because
 $$
 \x=\big(x_0,f^{-1}_{T(0)}(x_0),
  f^{-1}_{T(1)}\big(f^{-1}_{T(0)}(x_0)\big),f^{-1}_{T(2)}\big(f^{-1}_{T(1)}
  (f^{-1}_{T(0)}(x_0))\big),...\big).
 $$
\label{typeargument}
\noindent In other words,  identifying $\{0,1\}^{\N}$ with Cantor set $\CC$,
we have defined an injective map from $\X$ into $[0,1]\times \CC$,
but we must stress that this map is neither surjective nor
continuous, so  in order to get a precise representation of $\X$ we need to  develop
some new tools. We start with the following

\vspace{3mm}

\noindent {\it Definition} 1.
%\label{types}
Let $T$ be in $\{0,1\}^{\N}$ and let $\x=(x_0,x_1,...,x_n,...)\in \X$ be such that
$x_{n+1}=f^{-1}_{T(n)}(x_{n})$. Then we  say that $\x$ \emph{is of type} $T$, or
that $T$ \emph{is the type of} $\x$.
 We  denote by $\X_T\subset \X$ the set of all points  of  type $T$, and we say that
$T$ is \emph{admissible} for $(X,f)$ if $\X_T$ is not empty.

\vspace{3mm}

The coming  next lemma is straightforward. It says  that   the subsets $\X_T$ are
homoeomorphic to intervals. Loosely speaking these intervals are
 bricks with  help of which we shall build a homeomorphic image  of $\X$.

\begin{lem}[Bricks]\label{lem2}
Let $T\in \{0,1\}^{\N}$. Then the factor mapping
$$
\Phi(x_0,x_1,...)=x_0
$$

\noindent is a homeomorphism from $\X_T$  onto $\Phi(\X_T)\subset [0,1]$.
\end{lem}

As we already have 'bricks' we only need some 'cement' to start
our construction. A simple observation that for a positive number
$x$  close to zero,  inverse images $f_1^{-1}(x)$ and
$f_0^{-1}(x)$ are close to   each other, leads us to the following

\begin{lem}[Cement]\label{lem3}
Let $\x_0=(x_0,x_1,...)\in\X_T$ be such that $x_n=0$ for some $n\in\N$:
$$
\x_0=(f^n(0),f^{n-1}(0),...,f(0),0,...),
$$

\noindent and let
$T'\in \{0,1\}^{\N}$ be such that $T'(k)=T(k)$ for $k\neq n$, and $T'(n)=1\neq T(n)=0$.
Then   every neighborhood of $\x_0$ contains a subset
$$
\U=\{\x=(x_0',x_1',...)\in\X_{T'}: x_0'\in (f^n(0)- \varepsilon,f^n(0)+
\varepsilon)\}\subset \X_{T'}
$$

\noindent for some  $\varepsilon>0$, that is $\X_T\cup \X_{T'}$ is connected  topological
subspace of $\X$.
\end{lem}

\vspace{3mm}

\noindent {\it Proof.}
Take a neighborhood $\V$ of $\x_0$ of the form
$\V=\{(x_0',x_1',...)\in\X: |x_m'- x_m|<\varepsilon_0\}$
for some $\varepsilon_0 >0$ and $m\geq n$.
Such sets generate a base of neighborhoods of $\x_0$.
Since for $(x_0',x_1',...)\in\X$ we have $x_k=f(x_{k+1})$, $k\in\N$, and $f$ is continuous,
  there exists $\varepsilon_1>0$ such that
$$
\V_1=\{(x_0',x_1',...)\in\X:
|x_{n+1}' - x_{n+1}|< \varepsilon_1,\,\, x_{k+1}'=f^{-1}_{T(k)}(x_{k}'),
 \,\,\text{for}\,\,k > n\} \subset \V.
$$

\noindent As $\lim_{h\rightarrow 0}f^{-1}_1(h)=\rho=x_{n+1}$, there exists $\varepsilon_2>0$
such that
$f^{-1}_1(0,\varepsilon_2)\subset (x_{n+1}-\varepsilon_1,x_{n+1}+\varepsilon_1)$,
 that is $\V_{2}=\{(x_0',x_1',...)\in\X_{T'}: 0< x_{n}' < \varepsilon_2\} \subset \V_1$.
 By the same argument used to ensure the existence of $\varepsilon_1$  there exists
 $\varepsilon > 0$
  with the property that $\V_3=\{(x_0',x_1',...)\in\X_{T'}: |x_{0}' -x_{0}|< \varepsilon\}
  \subset \V_2$.
  Taking $\U=\V_3$ ends the  proof.

  \hfill$\blacksquare$

\noindent {\it Definition} 2.
%\label{thisdefinition}
Let $\x_0$, $\X_T$, $\X_{T'}$ be as in the preceding lemma. We
shall say that sets $\X_T$ and $\X_{T'}$ are \textit{connected in}
$\x_0$, and write $\X_T \sim \X_{T'}$ in $\x_0$.

\vspace{3mm}

If   $\X_T \sim \X_{T'}$ in $\x_0\in \X_T$, then $\Phi(\x_0)$ is an end point of the
interval $\Phi(\X_T)$.
In view of the preceding lemma we can 'glue'  the intervals  $\Phi(\X_T)$ and $\Phi(\X_T')$
into one,
and thus obtain a homeomorphic image of  $\X_T\cup \X_{T'}$. This is the fact we shall use
frequently and this is the cause that
we shall define  homeomorphisms of inverse limits in an inductive way.

This scheme  is  especially fruitful when dealing with the set
$$
\X_1=\{(x_0,x_1,x_2...)\in \X: \lim_{n\rightarrow \infty} x_n=1\}.
$$

\noindent  As $1$ is a repellent point, it follows  that $\X_1$  is open in $\X$ and
could be alternatively
  defined as the sum of all sets $\X_T$ for which almost all indices of $T$ are  $1$.
  We shall show that $\X_1$ is homeomorphic to a semi-open  interval in all cases.

\subsection{Cases one and two: $\rho\geq \rho_1$}\label{4}

In this section we describe the system $(\X,f)$ associated with  the $\F_{2^n}$ mapping
$f$ given by (\ref{f(x)}),
in case $\rho\geq \rho_1$, cf. Lemma\,1. First we define a homeomorphism from
$\X_1$  onto a semi-open
interval $(a,1]$, see Proposition\,1. Then  we show that in case
$\rho >\rho_1$, $\X$ is homeomorphic
to a closed interval, see Theorem\,2, and  in case $\rho=\rho_1$,
$\X$ is homeomorphic to
a subset of $\mathbb{R}^2$ pictured on Fig.\,3, see Theorem\,3.

Let us begin with some notations. We  denote by  $T_\infty$ an element of $\{0,1\}^{\N}$
with all entries equal to $0$,
 and  by $T_n$ an element of $\{0,1\}^\mathbb{N}$ where $T_n(k)=0$  for $k=0,...,n-1$,
 and $T_n(k)=1$ for $k\geqslant n$, that is:
$$
T_\infty=(0,0,0,...)\qquad  \textrm{ and }\qquad  T_n=(\underbrace{0,0,...,0}_{n},1,1,1...),
 \,\,\textrm{ for } n\in \N.
$$

\noindent Let us observe  that these are all the types admissible for $(X,f)$
(see Definition\,1).
Indeed, if $\x=(x_0,x_1,...)\in\X$ is such that  $x_{n+1}=f^{-1}_1(x_{n})$, for some
$n\in \N$,
then   $x_{k+1}=f^{-1}_1(x_{k})$ for all $k>n$, which follows immediately from the
fact that $\rho_1\leq\rho< x_{n+1}=f^{-1}_1(x_{n})$.
 Thus we have
$$\X=\bigcup_{n=0}^{\infty}\X_{T_n}\cup \X_{T_\infty}\qquad  \textrm{ and }\qquad
\X_1=\bigcup_{n=0}^{\infty}\X_{T_n}.$$

\noindent It is not hard to see that $\Phi(\X_{T_0})=(0,1]$ and $\Phi(\X_{T_n})$
is the image
of $[0,\rho_1)$  under $f_0^{n-1}$ for $n>0$.
 By Lemma\,2,  $\X_{T_n}$ is homeomorphic to the corresponding interval, that is,
$$\X_{T_n}\cong [f_0^{n-1}(0), f_0^n(0))\,\, \textrm{ for odd } n,\quad
 \X_{T_n}\cong (f_0^n(0), f_0^{n-1}(0)]\,\, \textrm{ for even } n>0.
$$

\noindent Using Lemma\,3 one can 'stick' those intervals together into one,
because  for $k\in \N$ we have
\begin{equation}\label{eq1}
 \X_{T_{2k}}\sim \X_{T_{2k
 +1}}\,\,\textrm{ in }\,\,(f^{2k}(0),f^{2k-1}(0),...,0,...)\in \X_{T_{2k+1}},
\end{equation}
\begin{equation}\label{eq2}
   \X_{T_{2k+1}}\sim \X_{T_{2k+2}}\,\,\textrm{ in }\,\, (f^{2k+1}(0),f^{2k}(0),...,0,...)
   \in \X_{T_{2k+2}},
\end{equation}

\noindent cf. Definition\,2.
 \begin{prop}\label{X_1_in_case_1_and_2}
Let $(X,f)$ be as in case 1 or 2 of the Lemma \ref{lemma1}  The mapping
\begin{equation}\label{firstformula}
\Psi(\x)=(-1)^n \Phi(\x)+2\sum_{k=1}^{n-1}(-1)^{k} f^{k}(0),\qquad\quad \text{for}\,\,\, \x
\in \X_{T_n},
\end{equation}

\noindent is a homeomorphism from $\X_1$ onto (maybe infinite) interval $(a,1]\subset
(-\infty,1]$, see Fig.\,2.
\end{prop}

\vspace{3mm}

\noindent {\it Proof.}
 The only item except the formula (\ref{firstformula}) left for us to prove is that sets
 $\X_{T_n}$ and $\X_{T_m}$
  where $m>n+1$ can be separated by open sets in $\X$. To see that take for instance:
  $\V_1=\{(x_0,...)\in \X_1: x_{n+1}>\rho_1\}$, and $\V_2=\{(x_0,...)\in \X_1:
  x_{n+1}<\rho_1\}$,
  then we have $\X_{T_n}\subset \V_1$, $\X_{T_m}\subset \V_2$ and $ \V_1\cap \V_2=\emptyset$.

  We use the mathematical induction method in order to construct the homeomorphism
  $\Psi:\X_1\rightarrow  (a,1]$, see Fig.\,2.

 Define $\Psi$ on $\X_{T_0}$ by $\Psi=\Phi$, and on $\X_{T_1}$ by $\Psi=-\Phi$, then
 $\Psi:\X_{T_0}\cup\X_{T_1}\rightarrow (-\rho_1,1]$ is a homeomorphism, see (\ref{eq1})
 and Lemma\,3.

Now, suppose  we have constructed a homeomorphism $\Psi$ from $\bigcup_{k=0}^{n-1} \X_{T_k}$
onto $(a_{n-1},1]$
 for $n > 1$, where $a_{n-1}=2\sum_{k=1}^{n-2} (-1)^kf_0^k(0)+(-1)^{n-1}f_0^{n-1}(0)$.
\Item{1)}
If $n$ is odd, then  put  $\Psi(\x)=-\Phi(\x)+a_{n-1}+f^{n-1}(0)$ for $\x\in \X_{T_n}$.
Then   $\Psi$ maps $\bigcup_{k=0}^{n} \X_{T_k}$ onto $( a_{n-1} -(f_0^n(0)-f_0^{n-1}(0)),1 ]=
 (2\sum_{k=1}^{n-1} (-1)^kf_0^k(0)+(-1)^{n}f_0^{n}(0),1]$
 and by (\ref{eq1})  it is a homeomorphism.  For $\x\in \X_{T_n}$, $\Psi(\x)$ is given by
 (\ref{firstformula}).
\Item{2)}
If $n$ is even, then  put  $\Psi(\x)=\Phi(\x)+a_{n-1}-f^{n-1}(0)$ for $\x\in \X_{T_n}$.
Then   $\Psi$ maps $\bigcup_{k=0}^{n} \X_{T_k}$ onto $(a_{n-1} +(f_0^n(0)-f_0^{n-1}(0)),1 ]=
(2\sum_{k=1}^{n-1} (-1)^kf_0^k(0)+(-1)^{n}f_0^{n}(0),1]$
 and by (\ref{eq2}) it  is a homeomorphism. For $\x\in \X_{T_n}$, $\Psi(\x)$ is given by
 (\ref{firstformula}).

In this manner  we obtain that  $\Psi$ is a homeomorphism from
$\X_1=\bigcup_{k=0}^{\infty} \X_{T_k}$
onto the interval $(a,1]\subset(-\infty,1]$ where $\displaystyle{a=2\sum_{n=0}^\infty
f^{2n}(\rho_1)-f^{2n+1}(\rho_1)}$,
and the formula (\ref{firstformula}) holds.

\begin{center}
\setlength{\unitlength}{0.8mm}
\begin{picture}(110,44)(-95,-5)
\thicklines
\qbezier(-35,11.5)(0,11.5)(36,11.5)\put(36,11.5){\circle*{0.8}}\put(35,7.5){\scriptsize{$1$}}
\put(3,11.5){\circle*{0.8}}\put(2.5,7.5){\scriptsize{$0$}}
\put(-16,11.5){\circle*{0.8}}\put(-18.5,7.5){\scriptsize{$-\rho_1$}}
\put(-35,11.5){\circle*{0.8}}\put(-37.5,7.5){\scriptsize{$-2\rho_1$}}
\put(-43,11.4){.}\put(-52,11.4){.}\put(-60,11.4){.}
\put(-70,11.5){\circle*{0.8}}\qbezier(-70,11.5)(-78,11.5)(-89,11.5) \put(-89,11.5){\circle*{0.8}}
\put(-91.5,10.5){\scriptsize{$\underbrace{\,\,\,\qquad\qquad\quad}_{|f^{n}(0)-f^{n+1}(0)|}$}}

\put(-117,11.5){\circle*{0.8}}\put(-117.5,8){\scriptsize{$a$}}
\thinlines
\qbezier[100](-35,11.5)(-70,11.5)(-115,11.5)
\put(16,32.5){$\X_{T_0}$}\put(18,30){\vector(0,-1){16}}\put(20,22){\scriptsize{$\Phi$}}
\put(-9,32.5){$\X_{T_1}$}\put(-7,30){\vector(0,-1){16}}\put(-5,22){\scriptsize{$-\Phi$}}
\put(-28,32.5){$\X_{T_2}$}\put(-26,30){\vector(0,-1){16}}\put(-24,22){\scriptsize{$\Phi- 2\rho_1$}}
\put(-81,32.5){$\X_{T_{n}}$}\put(-79,30){\vector(0,-1){16}}
\put(-77,22){\scriptsize{$($-$1)^{n}\Big(\Phi -f^{n-1}(0)\Big)$+$a_{n-1}$}}
\put(-46,-4){\small{Fig. 2.}}
\end{picture}
\\
\end{center}
Let us notice that if $\rho>\rho_1$, then $\omega_0$ is
 attractive and as $\omega_0\in (0,\rho_1)$, the sequence $\{\Phi(\X_{T_n})\}_{n>0}$
forms a decreasing family of intervals. We have $\lim_{n\rightarrow\infty}|\Phi(\X_{T_n})|=0$
where $|A|$ denotes
a diameter of a set $A\subset [0,1]$.  \\
On the other hand, if  $\rho=\rho_1$,  then $\Phi(\X_{T_n})=[0, \rho)$ for odd $n$, and
$\Phi(\X_{T_n})=(0, \rho]$ for even $n>0$,
 hence $|\Phi(\X_{T_n})|=\rho$ for all $n>0$ and as a consequence $\Psi(\X_1)=(-\infty,1]$.

\subsubsection{The case of attractive point $\omega_0$: $\rho>\rho_1$}

In order to get a complete homoeomorphic image of the inverse limit space $\X$ we need
to extend the homeomorphism from
Proposition \ref{X_1_in_case_1_and_2} onto the remaining  part $\X\setminus{\X_1}=
\X_{T_\infty}$ of $\X$. In case  $\rho>\rho_1$ it is
  very easy because  by attractiveness of $\omega_0$,   $\X_{T_\infty}$ is a singleton
  $\{\Tomega_0\}$ where $\Tomega_0$ is the  fixed point $(\omega_0,\omega_0,...)$ for $\f$.

\begin{thm}\label{atrractivepointtheorem}
Let $(X,f)$ be as in case one of Lemma\,1. Then up to
topological conjugacy  $\X$  is an interval $[0,1]$ and $\f$ is  a
monotonely increasing mapping with fixed points $0$ and $1$.
\end{thm}

\noindent {\it Proof.}
First, let us observe that if $\X= [0,1]$ then, since $\f$ is a homeomorphism
of $\X$ with two fixed points, it has the desired form.
Hence, it suffices to   show  that $\X$ is homeomorphic to a closed interval.

 As we have mentioned $\X\setminus\{\Tomega_0\}=\X_1$.  By Proposition
 \ref{X_1_in_case_1_and_2},
 $\X\setminus\{\Tomega_0\}$ is homeomorphic to $(a,1]$. Thus
defining $\Psi(\Tomega_0)=a$ we obtain a bijective mapping from $\X$ onto $[a,1]$,
and  to complete a task
we only need to investigate the  neighborhoods  of $\Tomega_0$.

For that purpose observe that sets $O_n(\Tomega_0)=\{(x_0,...)\in
\X: x_n\in (0,\rho)\}$, $n>0$, form  a base of neighborhoods of
$\Tomega_0$. Indeed, by  attractiveness of $\omega_0 \in(0,\rho)$,
$\bigcap_{n\in \N}f^n(0,\rho)=\{\omega_0\}$ and hence every
neighborhood of $\Tomega_0$ contains all $O_n(\Tomega_0)$ for $n$
greater than a certain $N$.
  Furthermore,  we
  have
$$
O_n(\Tomega_0)=\{\Tomega_0\}\cup \Big(\bigcup_{k=n}^\infty \X_{T_k}\Big)\setminus
\{(f^n(0),f^{n-1}(0),...,f(0),0,...)\},\qquad n>0,
$$

\noindent The form of the mapping $\Psi$ constructed in Proposition
\ref{X_1_in_case_1_and_2} implies that
 $$\Psi(\Big(\bigcup_{k=n}^\infty \X_{T_k}\Big)\setminus
 \{(f^n(0),f^{n-1}(0),...,f(0),0,...)\})=(a,b_n)\ {\rm for\ some}\ b_n\in(a,1]$$

\noindent and as $\Psi(\Tomega_0)=a$ we conclude that $\Phi$ is  a  homeomorphism
from $\X$ onto $[a,1]$.
Changing the scale on the real axis we can write $\X\cong [0,1]$.

\hfill$\blacksquare$

\subsubsection{The case of symmetry of the left hyperbola: $\rho=\rho_1$}

 Let us assume  now that $\rho=\rho_1$, see Lemma\,1.

 Firstly, let us observe that $f_0$ has the form
 $ f_0(x)=\displaystyle{\frac{\al x- \al\rho}{\gamma x -\al}},$
 $x\in [0,\rho],$ and if we conjugate $f$ with the homeomorphism
$$
h(x)=\left \{ \begin{array}{ll} \displaystyle{\frac{(c\rho+d) x}{c x+ d}} , & x\in [0,\rho],\\
x , & x\in (\rho,1],\end{array} \right.
$$

\noindent where  $c=\frac{1}{\rho}[1-(1-\gamma\frac{\rho}{\al})^{\frac{1}{2}}]$ and
$d=(1-\gamma\frac{\rho}{\al})^{\frac{1}{2}}$,
then the result mapping $g=h\circ f\circ h^{-1}$ acts as follows:
$g(x)=\left \{ \begin{array}{ll} \rho - x , & x\in [0,\rho],\\
 \frac{x-\rho}{1-\rho}, & x\in (\rho,1],\end{array}\right.$.
Hence, we may  assume  that $\al=-1$ and $\gamma=0$, that is
\begin{equation}
f(x)=\left \{ \begin{array}{ll} f_0(x)=\rho - x , & x\in [0,\rho],\\
f_1(x)=\displaystyle{\frac{x-\rho}{1-\rho}} , & x\in (\rho,1],\end{array} \right.
\end{equation}

\noindent Secondly, let us recall that in the case under consideration
$\X\setminus{\X_1}=\X_{T_\infty}$, and the  mapping $\Psi$ establishes a
homeomorphism between  $\X_1$ and  $(-\infty,1]$, see Proposition\,1.

Thirdly, let us notice that applying  Lemma\,2 to $\X_{T_\infty}$ we see that $\Phi$
is a
homeomorphism from $\X_{T_\infty}$ onto $[0,\rho]$. Hence defining
\begin{equation}\label{secondformula}
\Psi(\x)=(\Phi(\x),\infty), \qquad \text{for}\quad \x\in \X\setminus{\X_1},
\end{equation}

\noindent we conclude that $\Psi$ is a bijection from $\X$ onto $(-\infty,1]\cup
[0,\rho]\times\{\infty\}$,
and  $\Psi$ establishes  a homeomorphism from $\X\setminus{\X_1}$ onto the  copy
$[0,\rho]\times\{\infty\}$ of the interval $[0,\rho]$.

Finally, let us suppose that $(-\infty,1]\cup
[0,\rho]\times\{\infty\}$ is equipped with the topology carried over by
$\Psi$ from $\X$. Then $\Psi$ becomes a homeomorphism and the
system $(\X,\f)$ conjugated by $\Psi$ is described  in the coming
statement. By $[\cdot]$ we denote  here \emph{entier} function.

\begin{thm}\label{rho_1=rhotheorem}
Let $(X,f)$ be as in the case two of the  Lemma\,1. Then up to  topological
conjugacy
$\X$  consists of $\X_1=(-\infty,1]$ and a set $\X\setminus\X_1=[0,\rho]\times \{\infty\}$
naturally
homeomorphic to $[0,\rho]$  where  every neighborhood of point $(x,\infty)\in [0,\rho]\times
\{\infty\}$ contains,
for some $\varepsilon>0$ and $N\in\N$, all the intervals
\begin{equation}\label{neighbor0}
((-1)^{n} x - 2\Big[\frac{n}{2}\Big]\rho -\varepsilon,(-1)^{n} x - 2\Big[\frac{n}{2}\Big]
\rho +\varepsilon) \subset
(-\infty,1],\quad \text{for  } n > N.
\end{equation}

\noindent Hence $\X$ can be considered as a subset of $\mathbb{R}^2$
pictured on Fig.\,3. Mapping $\f$ acts on $\X$ as follows: It is  is
monotonely increasing on $(-\infty,1]$ where  $1$ is its fixed
point, namely
$$\f(x)=\left \{ \begin{array}{ll} \displaystyle{\frac{x-\rho}{1-\rho}} & x\in(\rho,1],\\
x-\rho , & x\in(\infty,\rho],\end{array} \right.
 $$

\noindent and  on $[0,\rho]\times \{\infty\}$, $\f$  is an involution, that is
$\f(x,\infty)=(\rho-x,\infty)$, for all
$(x,\infty)\in [0,\rho]\times \{\infty\}$.
\end{thm}

\vspace{3mm}

\noindent {\it Proof.}
Let $(-\infty,1]\cup[0,\rho]\times \{\infty\}$ be equipped with
the topology described above. We  show that $\Psi:\X\rightarrow
(-\infty,1]\cup [0,\rho]\times \{\infty\}$ is a homeomorphism.

We know that  $\Psi$ restricted to $\X_1$ is a homeomorphism onto
$(-\infty,1]$ and restricted to $\X\setminus\ X_1$ is a
homeomorphism  onto $[0,\rho]\times \{\infty\}$. Thus, since
$\X_1$ is open in $\X$, we only need to investigate neighborhoods
of points from $\X\setminus\ X_1$.

Let $\x_0=(x_0,f^{-1}_0(x_0),x_0,f^{-1}_0(x_0),...)$ be in  $\X\setminus\ X_1$.
It is obvious that sets
$
O_{n,\varepsilon}(\x_0)=\{(y_0,y_1...)\in \X: y_{2n}\in (x_0-\varepsilon,x_0+\varepsilon)\},
$
where  $n\in \N$, $\varepsilon >0$, form a base of neighborhoods of $\x_0$. For sufficiently
small $\varepsilon$ we obtain
$$
O_{n,\varepsilon}(\x_0)=\{(y_0,y_1...)\in \bigcup_{k=2n}^{\infty}\X_{T_k}\cup\X_{T_\infty}:
 y_{0}\in (x_0-\varepsilon,x_0+\varepsilon)\},
$$

\noindent and thus by (\ref{firstformula})  we have
$$
\Psi(O_{n,\varepsilon}(\x_0)\cap \X_1)=\bigcup_{k=2n}^{\infty}((-1)^{k} x_0 - 2\Big[
\frac{k}{2}\Big]\rho -\varepsilon,(-1)^{k}
x_0 - 2\Big[\frac{k}{2}\Big]\rho +\varepsilon).
$$

\noindent Hence $\Psi$ is a homeomorphism.

The form of the homeomorphism $\Psi\circ f\circ \Psi^{-1}$ can be now verified easily.
 We only mention that as $\rho_1=\rho$, the formula (\ref{firstformula}) simplifies
 to the following one
$$
\Psi(\x)=(-1)^n \Phi(\x) - 2[\frac{n}{2}]\,\rho,\ \ {\rm for}\ \ \x\in \X_{T_n}.$$

\hfill$\blacksquare$

\begin{center}
\setlength{\unitlength}{1mm}
\begin{picture}(110,65)(-85,-18)
\thinlines
\qbezier(-4,7.5)(9.5,18)(23,28.5)\put(23,28.5){\circle*{1}}\put(24,29){$1$}
\qbezier(-4,7.5)(-12,20)(-20,32.5)\put(-4,7.5){\circle*{1}}\put(-2,4.5){$0$}
\qbezier(-20,32.5)(-27,19.25)(-34,6)\put(-20,32.5){\circle*{1}}\put(-17,31){-$\rho$}
\qbezier(-34,6)(-40,19.75)(-46,33.5)\put(-34,6){\circle*{1}}\put(-31,4.5){-$2\rho$}
\qbezier(-46,33.5)(-51,19.5)(-56,5.5)\put(-46,33.5){\circle*{1}}\put(-43,32){-$3\rho$}
\qbezier(-56,5.5)(-60,20)(-64,34.5)\put(-56,5.5){\circle*{1}}\put(-53,4){-$4\rho$}
\qbezier(-64,34.5)(-67,20.25)(-70,5)\put(-64,34.5){\circle*{1}}
\qbezier(-70,5)(-72,20)(-74,35)
\qbezier(-76,5)(-75,20)(-74,35)
\qbezier(-76,5)(-76.5,20)(-77,35)
\qbezier[150](-78,5)(-77.5,20)(-77,35)
\qbezier[100](-78,5)(-78.5,20)(-79,35)
\qbezier[75](-79,5)(-79,20)(-79,35)

\qbezier(-80,5)(-80,20)(-80,35)
\put(-80,5){\circle*{1}}\put(-86,1){$(0,\infty)$}
\put(-80,35){\circle*{1}}\put(-86,37){$(\rho,\infty)$}
\end{picture}

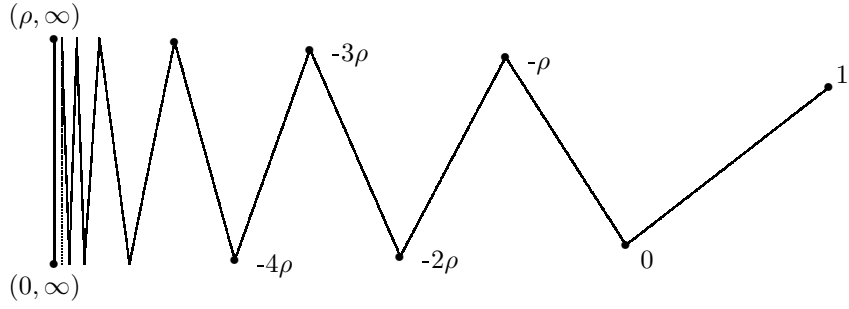
{\small Fig.\,3. The inverse limit
associated with $([0,1],f)$ in case $f$ is $\F_{2^1}$ mapping
possessing an interval  of periodic
points.}

\end{center}

\subsection{Case three: $\rho < \rho_1$}\label{5}

In this section we describe the inverse limit space $\X$ and the induced homeomorphism $\f$
in case $\rho < \rho_1$, cf. Lemma\,1. First we find all types admissible
for $(X,f)$, cf. Definition\,1.
Then we construct a homeomorphism from $\X_1$  onto a semi-open  interval $(a,1]$,
see Proposition\,2.
 Afterwards   we describe the subset $\X_{\omega_0}$ of $\X$ consisting of elements
 with limit point $\omega_0$.
 Namely, we show that   $\X_{\omega_0}$ is homeomorphic to an open interval, see
 Proposition\,3.
 Finally we gather  these facts within Theorems\,4 and
5 and obtain a description
 of $(\X,\f)$
     in cases 3a) and 3b) of Lemma\,1.

Throughout this section we assume that $\rho_1>\rho$. Then there
are two repellent fixed points $\omega_0$, $1$, and two  periodic
points $w_1$, $w_2$ of period two, see Fig.\,4, recall  Lemma\,1.
In case a) the points $w_1$ and  $w_2$ are
attractive, and in case b) points from sets
$[0,f(\rho_1)]\setminus\{w_1\}$ and
$[\rho,\rho_1]\setminus\{w_2\}$  are of period four.

Let us find first  all types admissible for $(X,f)$, cf. Definition\,1.

To this end, let us  observe that since $f^2(\rho,\rho_1)\subset (\rho,\rho_1)$,
we have $\rho<f^2(\rho_1)$, and it follows  that
$f^{-1}_1(f^{-1}_0(\rho))>\rho_1$, and all the more $f^{-2}_1(\varepsilon)>\rho_1$,
for $\varepsilon>0$.
Thereby, if $T(n)=T(n+1)=0$ and $T(n+2)=1$, for some $n\in\N$, then we have $T(n+m)=1$
for all $m>1$, because
$f^{-1}_1(f^{-2}_0([0,\rho_1]))=f^{-1}_1([f^{-1}_0(\rho),\rho])\subset (\rho_1,1)$.
If $T(n)=T(n+1)=1$, for some $n\in\N$,
then $T(n+m)=1$ for all $m\in \N$ because $f^{-2}_1((0,1])\subset (\rho_1,1]$.

Concluding: If $T$ is admissible for $(X,f)$, then we have only the following  possibilities:

1) $T$  consists of $0$ and $1$ appearing alternately,

2) $T$  consists of $0$ and $1$ appearing alternately and
starting from  some $n$, $T$ is constant,

3) $T$  consists of $0$
and $1$ appearing alternately and for some $n$ we have
$T(n)=T(n+1)=0$ and $T(m)=1$ for $m>n+1$.

\begin{center}
\setlength{\unitlength}{1.1mm}
\begin{picture}(50,56)(0,-3)
\thinlines
\put(0,0){\framebox(50,50){}}
\qbezier[200](0,0)(25,25)(50,50)
\thicklines
\qbezier(0,33.5)(29.5,30)(30,0)
\qbezier(30,0)(40,25)(50,50)
\thinlines
\put(-1,-4){$0$}
\put(28,-3){$\rho$}
\put(49,-4){$1$}
\put(-5,33){$\rho_1$}\put(-0.5,32.5){-}\put(-4,28){$\rho$} \put(-0.5,28){-}
\put(-11.5,8){$f(\rho_1)$} \put(-0.5,8){-}

\qbezier[80](22.7,0)(22.7,11.25)(22.7,22.5)\put(21.5,-3){$\omega_0$}
\put(-3,49){$1$}

\qbezier[110](5.2,0)(5.2,16)(5.2,32.7)\put(4,-3){$w_1$}
\qbezier[20](32.2,0)(32.2,2.5)(32.2,5)\put(31.5,-3){$w_2$}
\qbezier[50](5.2,32.5)(13.5,32.5)(32,32.5)
\qbezier[50](32.1,5)(32.1,13.5)(32.1,32)
\qbezier[50](5.2,5.2)(13.5,5.2)(32,5.2)
\end{picture}

\vspace{3mm}

{\small Fig.\,4.}
\end{center}
\label{discussion}

From the above it follows  that  elements of $\X_1$ are either of type
from item 3) or from item 2) where almost all indices of $T$ are
equal $1$, cf. (\ref{X_1_equation}). Similarly the points from
$\X_{\omega_0}$, that is   elements with limit point $\omega_0$,
are of the type from the item 2) where starting from some place all indices of $T$
are equal $0$, cf. (\ref{X_omega0_i_takietam}). There are only two
types $T^1_\infty=(1,0,1,0,1,...)$ and $T^0_\infty=(0,1,0,1,0...)$
which  fall into  the item 1). It is not hard to see that in case $w_1$
and  $w_2$ are  attractive, sets $\X_{T^0_\infty}$,
$\X_{T^1_\infty}$ reduce to singletons, and in the case
$[0,f(\rho_1)]$ and $[\rho,\rho_1]$  are intervals of periodic
points. Sets $\X_{T^1_\infty}$, $\X_{T^0_\infty}$ are homeomorphic
to those intervals. Namely, we have either
$$
\Phi(\X_{T^0_\infty})=\{w_2\},\,\, \Phi(\X_{T^1_\infty})=\{w_1\}\,\textrm{ or }\,
\Phi(\X_{T^0_\infty})=
[\rho,\rho_1],\,\, \Phi(\X_{T^1_\infty})=[0,f(\rho_1)].
$$

\subsubsection{The set $\X_1$ of elements with limit point $1$}

In order to  investigate the subset $\X_1$ of $\X$ we provide the following  notation
$$
T_{n,k}^{(1)}=(\underbrace{\overbrace{1,0,1,...0,1,0}^{2k},0,...0,0}_{n},1,1,1...),
\qquad \textrm{ for } k=1,...,\Big[\frac{n}{2}\Big]
$$
$$
T_{n,k}^{(0)}=(\underbrace{\overbrace{0,1,0,...1,0,1}^{2k},0,...0,0}_{n},1,1,1...),\qquad
\textrm{ for } k=0,...,\Big[\frac{n-1}{2}\Big]
$$

\noindent where $n\in\N$. Within  the notation from the previous section we have
$T_n=T_{n,0}^{(0)}$ for $n>0$.
 In view of the  argument  from  the present  section we have
\begin{equation}\label{X_1_equation}
\X_1=\bigcup_{n\in\N}\,\bigcup_{i=0}^1\bigcup_{k=i}^{[\frac{n-1+i}{2}]}\X_{T_{n,k}^{(i)}}.
\end{equation}
Simple geometrical considerations lead to conclusion that $\Phi(\X_{T_0})=(0,1]$,
$\Phi(\X_{T_{1}})=[0,\rho_1)$,
and  for $n>1$, $i=0,1$,
 $$
 \Phi(\X_{T_{n}})=[0,\rho_1],
 $$
 $$
  \Phi(\X_{T_{2n,n}^{(1)}})=f^{2n-1}((\rho,\rho_1)), \qquad \Phi(\X_{T_{2n+1,n}^{(0)}})=
  f^{2n}((\rho,\rho_1)),
 $$
 $$
 \Phi(\X_{T_{n,k}^{(i)}})=f^{2k-i}((\rho,\rho_1]),\qquad \text{for}\,\,k=i,...,\Big[
 \frac{n-1+i}{2}\Big]-1.
 $$
By the Lemma\,2 sets $\X_{T_{n,k}^{(i)}}$ are homeomorphic to  corresponding
intervals.

Let us now introduce an order in the family of sets
$\X_{T_{n,k}^{(i)}}$, or equivalently in the family of types
$T_{n,k}^{(i)}$, such that any two neighboring with respect to
this order sets  are connected in the sense of the Lemma\,3.
For that purpose, we arrange them into clusters of $4n+2$
elements.

 The zero cluster is: $T_0$, $T_{1}$,  (recall that $T_n=T_{n,0}^{(0)}$), and for $n>0$
 the $n$-th cluster is
$$
T_{2n},\,\,T_{2n,1}^{(1)},\,\,T_{2n,2}^{(1)},\,...,\,\,
T_{2n,n}^{(1)},\,\,T_{2n+1,n}^{(1)},\,\,T_{2n+1,n-1}^{(1)},\,...,\,\,T_{2n+1,1}^{(1)},
$$
$$
T_{2n+1},\,\,T_{2n+1,1}^{(0)},\,\,T_{2n+1,2}^{(0)},\,...,\,\,T_{2n+1,n}^{(0)},\,\,
T_{2n+2,n}^{(0)},\,\,T_{2n+2,n-1}^{(0)},\,...,\,\,T_{2n+2,1}^{(0)}.
$$
Then every two types neighboring in the above cluster differ from each other with exactly
 one index, and they are connected in the sense of Lemma\,3, cf. Definition\,2.
One can check that
$$
\X_{T_{2n}} \sim \X_{T_{2n,1}^{(1)}} \textrm{ in  } (0,...)\in \X_{T_{2n}};
$$
$$
\X_{T_{2n,k}^{(i)}} \sim \X_{T_{2n,k+1}^{(i)}} \textrm{ in } (f^{2k-i}(\rho_1),...,
\rho_1,0,...)\in \X_{T_{2n,k}^{(i)}},
 \textrm{ for } k=i,...,\Big[\frac{n-1+i}{2}\Big]-1;
$$
$$
\X_{T_{2n,n}^{(1)}} \sim \X_{T_{2n+1,n}^{(1)}}  \textrm{ in } (f^{2n-1}(\rho_1),...,
\rho_1,0,...)\in \X_{T_{2n+1,n}^{(1)}};
$$
$$
\X_{T_{2n+1,1}^{(1)}} \sim \X_{T_{2n+1}} \textrm{ in } (\rho_1,0,...)\in \X_{T_{2n+1}}
\textrm{ and }
 \X_{T_{2n+1}}\sim \X_{T_{2n+1,1}^{(0)}}  \textrm{ in }   (0,...) \in \X_{T_{2n+1}}.
$$
$$
\X_{T_{2n+1,n}^{(0)}} \sim \X_{T_{2n+2,n}^{(0)}}  \textrm{ in } (f^{2n}(\rho_1),...,
\rho_1,0,...)\in \X_{T_{2n+2,n}^{(0)}};
$$
Moreover the first set  in the $n$-th cluster is connected with the last one in the
$(n-1)$-th cluster,
that is:  $\X_{T_{2n,1}^{(0)}} \sim \X_{T_{2n}}$ in $(\rho_1,0,...)\in \X_{T_{2n}}$.

Thus we have arranged sets $\X_{T_{n,k}^{(i)}}$ in the following way:
$$
\X_{T_0},\,\,\X_{T_{1}},\,\,\X_{T_{2}},\,\,\X_{T_{2,1}^{(1)}},\,\,\X_{T_{3,1}^{(1)}},\,\,
\X_{T_{3}},\,\,\X_{T_{3,1}^{(0)}},\,\,\X_{T_{4,1}^{(0)}},\,\,\X_{T_{4}},\,\,
\X_{T_{4,1}^{(1)}},\,\,\X_{T_{4,2}^{(1)}},$$
$$
\X_{T_{5,2}^{(1)}},\,\,\X_{T_{5,1}^{(1)}},\,\,\X_{T_{5}},\,\,\X_{T_{5,1}^{(0)}},\,\,
\X_{T_{5,2}^{(0)}},\,\,\X_{T_{6,2}^{(0)}},\,\,\X_{T_{6,1}^{(0)}},\,\,\X_{T_6},\,\,...
$$

\noindent Let us denote
 $$ A_0=\X_{T_0},\,\,\, A_1=\X_{T_{1}},\,\,\, A_2=\X_{T_{2}},\,\,\, A_3=\X_{T_{2,1}^{(1)}},
 \,\,\,
 A_4=\X_{T_{3,1}^{(1)}},\,\,\, A_5=\X_{T_{3}},\,\,\, ...\,.
$$

\noindent Thereby we have $\X_1=\bigcup_{n\in \N}A_n$.  We denote by $|\Phi(A_n)|$ the length of
interval $|\Phi(A_n)|$,
 and by $dist(\Phi(A_n),x)$ the distance between set $\Phi(A_n)$ and point $x\in\R$.

\begin{prop}\label{X_1_in_case_3}
Let $(X,f)$ be as in the case 3 of the Lemma\,1. The mapping
\begin{equation}\label{sickhomeomorphism}
\Psi(\x)=(-1)^n\Big(\Phi(\x)-dist(\Phi(A_n),0)\Big)- \sum_{k=1}^{2[\frac{n}{2}]}|\Phi(A_k)|,
\quad \textrm{for }\x \in A_n,
\end{equation}
is  a homeomorphism from $\X_1$ onto  an interval $(a,1]\subset(-\infty,1]$.
\end{prop}

\vspace{3mm}

\noindent {\it Proof.}
First observe that if $|i-j| > 1$, then sets $A_i=\X_T$ and $A_j=\X_{T'}$ possess disjoint
open neighborhoods $\V_1$ and
$\V_2$ in $\X_{1}$, that is  $A_i\subset \V_1$, $A_{j}\subset \V_{2}$, and $\V_1\cap\V_2=
\emptyset$.
Indeed, as $|i-j| > 1$,  sequence $T=T_{m,k}^{i}$ differs from $T'=T_{m',k'}^{i'}$ at least
with two entries,
 and in view of the discussion from page \pageref{discussion} the only  cases we need to
 consider are:
\begin{Item}{1)}
If for some $n$, $T(n)=T(n+1)=2$ and $T'(n)=T'(n+1)=1$, then we take $\V_1=\{(x_0,...)\in
\X_1: x_{n+2}>\rho_1\}$
and $\V_2=\{(x_0,...)\in \X_1: x_{n+2}<\rho_1\}$.
\end{Item}
\begin{Item}{2)}
If for some $n$, $T(n)=T(n+2)=2$, $T(n+1)=T(n+3)=1$, $T'(n)=T'(n+1)=T'(n+2)=T'(n+3)=1$
then we can take
$\V_1=\{(x_0,...)\in \X_1: x_{n+3}>\rho\}$
and $\V_2=\{(x_0,...)\in \X_1: x_{n+3}<\rho\}$.
\end{Item}
 Now, using Lemma \ref{lem3}, we can construct a homeomorphism $\Psi: \X_1\rightarrow
 (-a_0,1]$.

The recipe: put $\Psi(\x)=\Phi(\x)$ for $\x\in\X_{T_0}$,  $\Psi(\x)=-\Phi(\x)$ for
$\x\in\X_{T_1}$,
and $\Psi(\x)=\Phi(\x) -2\rho_1$ for $\x\in\X_{T_2}$ then clearly $\Psi:\X_{T_0}
\cup\X_{T_1}\cup\X_{T_2}
\rightarrow [-2\rho_1,1]$ is a homeomorphism.

 Now let $n>2$ and suppose  we have defined  homeomorphic mapping $\Psi$ from
 $\bigcup_{k=0}^{n-1} A_k$
 onto an interval $[a_{n-1},1]$ or $(a_{n-1},1]$ (depending on $n$), where $a_{n-1} =
 -\sum_{k=1}^{n-1}|\Phi(A_k)|$:
 \begin{Item}{1)}
If $n$ is even (this is the case iff $A_n=\X_{T_{m,k}}^{(i)}$ where $m+k$ is even) then
for $\x\in A_n$
we  define $\Psi(\x)=\Phi(\x)-|\Phi(A_n)|-dist(A_n,0)$.
\end{Item}
\begin{Item}{2)}
If $n$ is odd (this is the case iff $A_n=\X_{T_{m,k}}^{(i)}$ where $m+k$ is odd), then for
$\x\in A_n$
 we  define $\Psi(\x)= -\Phi(\x)+dist(A_n,0)$.
\end{Item}
It can be easily checked that  we obtain in this manner a homeomorphism from
$\bigcup_{k=0}^{n} A_k$ onto
 $[-b,1]$ or $(-b,1]$ (depending on $n$) where $b=a+|\Phi(A_n)|$. Hence, the hypotheses
 follows.

\hfill$\blacksquare$

\noindent {\it Remark} 1.
%1
The homeomorphism from  the preceding proposition is constructed in the same fashion  as the
one from the
Proposition\,1. To see that note that in the case considered in the
latter proposition,
cf. (\ref{firstformula}), we have
$$
(-1)^{n+1}dist(\Phi(\X_{T_n}),0)- \sum_{k=1}^{2[\frac{n}{2}]}|\phi(\X_{T_k})|
= 2\sum_{k=1}^{n-1}(-1)^{k} f^{k}(0).
$$

\vspace{3mm}

In order to better understand  the action of homeomorphism (\ref{sickhomeomorphism}) let us
denote
by $d_n$, $n>0$,  the length  of the interval $f^n((\rho,\rho_1])$,  and by $d_0$ the length
of $[0,\rho_1]$,  that is
$$d_n=|f^{n}(\rho)-f^n(\rho_1)|,\qquad \textrm{ for } n>0,\qquad \textrm{ and }
\qquad d_0=\rho_1.$$

\noindent We recall that we split the sequence $\{A_k\}_{k\in\N}$ into
clusters of $4n+2$ elements:
$$A_0,\,A_1;\,\,\,A_2,\,...,A_{7};\,\,\, A_8,\,...,\,A_{17};\,\,\, ...\,\,\, ;\,\,\, A_{2n^2},
\,...,\,A_{2n^2+4n+ 1};\,\,\, ...\, .$$
Then for the $n$-th such cluster, $n>0$, we have
$$
|\Phi(A_{2n^2})|=|\Phi(A_{2n^2+2n +1})|=d_0,
$$
$$
|\Phi(A_{2n^2+k})|=|\Phi(A_{2n^2+2n -k +1})|=|\Phi(A_{2n^2+2n +k +1})|=|\Phi(A_{2n^2+4n -k +2})|
=d_k,
$$

\noindent where $k=1,...,n$ and  the picture is as follows

\begin{center}
\setlength{\unitlength}{0.8mm}
\begin{picture}(160,48)(-120,-7)
\thicklines
\qbezier(-45,11.5)(-20,11.5)(36,11.5)
\qbezier(-55,11.5)(-45,11.5)(-13,11.5)
\qbezier(-55,11.5)(-65,11.5)(-73,11.5)
\qbezier(-97,11.5)(-80,11.5)(-73,11.5)
\qbezier(5,11.5)(-4,11.5)(-13,11.5)
\put(-47,11.5){\circle*{0.8}}
\put(-31,11.5){\circle*{0.8}}\put(-45.5,10){$\underbrace{\qquad \quad\, }_{d_1}$}
\put(-13,11.5){\circle*{0.8}}\put(-29.5,10){$\underbrace{\qquad \quad\,\,\, }_{d_0}$}
\put(-63,11.5){\circle*{0.8}}\put(-62,10){$\underbrace{\qquad \quad \,}_{d_1}$}
\put(-81,11.5){\circle*{0.8}}\put(-80,10){$\underbrace{\qquad \quad \,\,\,}_{d_0}$}
\put(-97,11.5){\circle*{0.8}}\put(-96,10){$\underbrace{\qquad \quad \,}_{d_1}$}

\put(-117,11.5){\circle*{0.8}}\put(-117.5,13){\scriptsize{$a$}}

\put(5,11.5){\circle*{0.8}}\put(-11,10){$\underbrace{\qquad \quad \,\,\,}_{d_0}$}
\put(36,11.5){\circle*{0.8}}

\put(5,13){\scriptsize{$0$}}
\put(35,13){\scriptsize{$1$}}

\thinlines

\qbezier[100](-45,11.5)(-70,11.5)(-115,11.5)

\put(20,32.5){$\X_{T_0}$}\put(22,30){\vector(0,-1){16}}\put(24,22){\scriptsize{$\Phi$}}

\put(-41,32.5){$\X_{T_{2,1}^{(1)}}$}\put(-39,30){\vector(0,-1){16}}\put(-37,22){\scriptsize{-$\Phi$-$2\rho_1$}}
\put(-23.5,32.5){$\X_{T_2}$}\put(-22,30){\vector(0,-1){16}}\put(-21,22){\scriptsize{$\Phi$-$2\rho_1$}}
\put(-6,32.5){$\X_{T_1}$}\put(-4,30){\vector(0,-1){16}}\put(-3,22){\scriptsize{-$\Phi$}}
\put(-56,32.5){$\X_{T_{3,1}^{(1)}}$}\put(-54,30){\vector(0,-1){16}}\put(-62,25){\scriptsize{$\Phi$-$2\rho_1$-$2f(\rho_1)$}}
\put(-74,32.5){$\X_{T_3}$}\put(-72,30){\vector(0,-1){16}}\put(-80,19){\scriptsize{-$\Phi$-$2\rho_1$-$2f(\rho_1)$}}

\put(-91,32.5){$\X_{T_{3,1}^{(0)}}$}\put(-89,30){\vector(0,-1){16}}
\put(-101,25){\scriptsize{$\Phi$-$3\rho_1$-$2f(\rho_1)$-$f^2(\rho_1)$}}
\put(-52,-5){\small{Fig. 5.}}
\end{picture}
\\
\end{center}

\subsubsection{The set $\X_{\omega_0}$ of elements with limit point $\omega_0$}

In this subsection we  examine the set $\X_{\omega_0}$ consisting of elements  $\x=(x_0,...)$
 for which $\omega_0$ is an accumulation point. Namely,
$$
\X_{\omega_0}=\{\x\in\X:\lim_{n\rightarrow \infty}x_n=\omega_0\}.
$$
In view of the earlier discussion (see page \pageref{typeargument}),
$\X_{\omega_0}$ is the sum of sets of the form $\X_T$ where almost
all indices of $T$ are  equal $1$.  More precisely, if we denote
$$
T_k^{1}=(\underbrace{1,0,1,...0,1,0}_{2k},0,0,0,...),\qquad T_{k}^{0}=
(\underbrace{0,1,0,...1,0,1}_{2k},0,0,0...),
$$
(we recall that $T_\infty=(0,0,0,....)$, $T^1_\infty=(1,0,1,0,1,...)$, $T^0_\infty=
(0,1,0,1,0...)$), then  we have
$$
\X\setminus \X_1=\bigcup_{i=0}^1\bigcup_{k=1}^\infty\X_{T^i_k}\cup \X_{T^i_\infty}
\cup\X_{T_\infty},
$$
and obviously
\begin{equation}\label{X_omega0_i_takietam}
\X_{\omega_0}=\bigcup_{i=0}^1\bigcup_{k=1}^\infty\X_{T^i_k}\cup\X_{T_\infty}, \qquad
\textrm{ and }\qquad \X\setminus
\X_1=\X_{\omega_0}\cup \X_{T^0_\infty}\cup \X_{T^1_\infty}.
\end{equation}
As our aim at the moment is a description of $\X_{\omega_0}$ let us
notice that
$$
\Phi(\X_{T_\infty})=[0,\rho_1],\,\quad
\Phi(\X_{T_k^i})= f^{2k  -i}(\rho,\rho_1].
$$

\noindent Hence $|\Phi(\X_{T_\infty})|=d_0$ and $|\Phi(\X_{T_k^i})|=d_{2k  -i}$
(see the previous subsection), and we have the following

\begin{prop}\label{X_omega_in_case_3}
Let $(X,f)$ be as in the case 3 of the Lemma\,1.  The mapping $\Theta$
defined as follows
$$
\Theta(\x)= \Phi(\x), \qquad\quad \text{for}\,\,\, \x\in \X_{T_\infty},
$$
\begin{equation}\label{thirdformula}
\Theta(\x)=(-1)^n \Phi(\x)+2\sum_{k=i-1}^{n-1}(-1)^k f^{2k-i}(\rho_1) ,\qquad\quad
\text{for}\,\,\, \x\in \X_{T_n^i},
\end{equation}

\noindent is a homeomorphism from $\X_{\omega_0}$ onto an interval $(a_\infty,b_\infty)\subset
(-\infty,+\infty)$.\\
Moreover, the conjugated mapping $\Theta\circ  \f \circ
\Theta^{-1}$ is a decreasing  homeomorphism of
$(a_\infty,b_\infty)$ such that  interval of length $d_0$, marked
on Fig.\,6, is mapped onto itself and  interval of length $d_1$;
  interval of length $d_{2n}$, $n>0$,  is mapped onto  interval of length $d_{2n+1}$,
  and  interval of length $d_{2n-1}$  is mapped
  onto  interval of length $d_{2n}$. The unique fixed point of this mapping is $\omega_0$.
\end{prop}

\begin{center}
\setlength{\unitlength}{0.8mm}
\begin{picture}(160,45)(-120,-5)
\thicklines

\qbezier(-55,11.5)(-45,11.5)(-13,11.5)
\qbezier(-55,11.5)(-65,11.5)(-73,11.5)
\qbezier(-91,11.5)(-80,11.5)(-73,11.5)
\qbezier(5,11.5)(-4,11.5)(-13,11.5)
\put(-55,11.5){\circle*{0.8}}
\put(-31,11.5){\circle*{0.8}}\put(-54,10){$\underbrace{\qquad\qquad \,\,\,\,\, }_{d_0}$}
\put(-13,11.5){\circle*{0.8}}\put(-30,10){$\underbrace{\qquad \quad \,\,\,\,}_{d_2}$}
\put(-73,11.5){\circle*{0.8}}\put(-72,10){$\underbrace{\qquad \quad \,\,\,\,}_{d_1}$}
\put(-91,11.5){\circle*{0.8}}\put(-90,10){$\underbrace{\qquad \quad \,\,\,\,}_{d_3}$}
\put(5,11.5){\circle*{0.8}}\put(-12,10){$\underbrace{\qquad \quad \,\,\,\,}_{d_4}$}

\put(-117,11.5){\circle*{0.8}}\put(-117.5,13){\scriptsize{$a_\infty$}}
\put(39,11.5){\circle*{0.8}}\put(38.5,13){\scriptsize{$b_\infty$}}
\put(-55.7,13){\scriptsize{$0$}}
\put(-32,13){\scriptsize{$\rho_1$}}

\thinlines
\qbezier[100](-45,11.5)(-20,11.5)(36,11.5)
\qbezier[100](-45,11.5)(-70,11.5)(-115,11.5)
\put(-45,32.5){$\X_{T_\infty}$}\put(-43,30){\vector(0,-1){16}}\put(-41,22){\scriptsize{$\Phi$}}
\put(-23.5,32.5){$\X_{T_1^0}$}\put(-22,30){\vector(0,-1){16}}\put(-20,22){\scriptsize{-$\Phi$+$2\rho_1$}}
\put(-6,32.5){$\X_{T_2^0}$}\put(-4,30){\vector(0,-1){16}}\put(-2,22){\scriptsize{$\Phi$+$2\rho_1$-2$f^2(\rho_1)$}}
\put(-66,32.5){$\X_{T_1^1}$}\put(-64,30){\vector(0,-1){16}}\put(-62,22){\scriptsize{-$\Phi$}}
\put(-84,32.5){$\X_{T_2^1}$}\put(-82,30){\vector(0,-1){16}}\put(-80,22){\scriptsize{$\Phi$-$2f(\rho_1)$}}

\put(-52,-5){\small{Fig. 6.}}
\end{picture}
\\
\end{center}

\noindent {\it Proof.}
As $[0,f(\rho_1)]\cap [\rho,\rho_1]=\emptyset$ we conclude that
for all $k,l\in \N$ sets $\X_{T_k^0}$ and $\X_{T_l^1}$ can be
separated by open sets in $\X$. Now, we show that sets
$\X_{T_k^i}$ and $\X_{T_l^i}$, $i=0,1$, such that $|k-l|>1$ also
possess disjoint open neighborhoods.

We may assume that $k>l$. The condition $k-l>1$ implies that there exists $n\in N$
such that $
T_k^i(n)=T_k^i(n+2)=2$, $T_k^i(n+1)=1,$ and $ T_l^i(n)=T_l^i(n+1)=T_l^i(n+2)=1.
$
Then it is easy to see that
$$\X_{T_k^i}\subset \{(x_0,...)\in \X:x_n+3\in  (\rho,\rho_1]\},
$$
$$
\X_{T_l^i}\subset \{(x_0,...)\in \X:x_n+3\in [0,f^{-2}_1(\rho)]\}.
$$
Since $\omega_0$ is repellent $f^{-2}_1(\rho)< \rho$. Thus taking $q$ such that
$f^{-2}_1(\rho)<q< \rho$  and setting
$$
\V_1=\{(x_0,...)\in \X:x_n+3 > q \}, \qquad \V_2=\{(x_0,...)\in \X:x_n+3 < q \}
$$
we prove our claim.

In the same manner one can prove that $\X_{T_\infty}$ and
$\X_{T_k^i}$ for  $k>1$  can  also be separated by open sets.

Applying the Lemma\,3 we obtain that the remaining pairs of
sets, that is: $\X_{T_\infty}$, $\X_{T_1^i}$, and  $\X_{T_k^i}$,
$\X_{T_{k+1}^i}$,  cannot be separated, because
$$
(0,...)\in \X_{T_\infty},\qquad \textrm{and}\qquad(f^{2k-i}(\rho_1),...,\rho_1,0,...)
\in \X_{T_k^i}.
$$
Now, similarly as in proofs of Propositions 2
and 3,  in view of Lemmas\,2 and 3  one easily
checks that
 formula (\ref{thirdformula}) defines  a homeomorphism, cf. Fig.6.
\\
It is obvious that homeomorphism $\f$ preserves $\X_{\omega_0}$. Furthermore we  have
$$
\f(\X_{T_\infty})=\X_{T_\infty}\cup \X_{T_{1}^{1}},\,\,\,\,\,\,\,
\f(\X_{T_n^{0}})=\X_{T_{n+1}^{1}},\,\,\,\,\,\,\, \f(\X_{T_n^{1}})=\X_{T_{n}^{0}},
$$
and $\f(\omega_0,\omega_0,...)=(\omega_0,\omega_0,...)\in \X_{T_\infty}$. From this the
hypotheses follows.

\hfill$\blacksquare$

\subsubsection{The case of attractive cycle of period two}

For each $n>0$,  the point $w_1$ lies in $f^{2n-1}((\rho,\rho_1])$
and the point $w_2$ lies in $f^{2n}((\rho,\rho_1])$. Hence, in
case $w_1$ and $w_2$ form an attractive cycle we have
$$
d_1> d_3 > d_5 > ...\, ,\qquad  d_2 > d_4 > d_6 > ...\,, \qquad \text{and}\qquad
\lim_{n\rightarrow \infty} d_n=0.
$$
It follows that $\X_{T^i_\infty}=\{\w_{2-i}\}$, $i=0,1$, where
$\w_{1}=(w_1,w_2,w_1,...)$ and $\w_{2}=(w_2,w_1,w_2,...)$. In
other words $\X\setminus \X_1=\X_{\omega_0}\cup\{\w_1,\w_2\}$ and
thus we can easily prolong the homeomorphism from the Proposition\,3
 to the set  $\X\setminus \X_1$.

\begin{prop}\label{X minus X_1 in case_3a)}
Let $(X,f)$ be as in the case 3a) of the Lemma\,1, and let $\Theta$
be the mapping from Proposition\,3.
 If we define
$$
\Theta(\w_1)=a_\infty , \qquad \textrm{ and }\quad \,\,\, \Theta(\w_2)=b_\infty,
$$
 then $\Theta$ is a homeomorphism from $\X\setminus \X_1$ onto the  closed
 interval $[a_\infty,b_\infty]\subset[-\infty,+\infty]$.
\end{prop}

\vspace{3mm}

\noindent {\it Proof.}
We shall  consider only the point $\w_1$, an argumentation concerning the point
$\w_2$ is similar.\\
It is clear that sets $O_n(\w_1)=\{(x_0,...)\in \X\setminus \X_1: x_{2n-1}\in(\rho,\rho_1)\}$,
$n>0$, form
 a base of neighborhoods of $\w_1$ in $\X\setminus \X_1$. Furthermore,  we
  have
$$
O_n(\w_1)=\{\w_1\}\cup \Big(\bigcup_{k=n}^\infty \X_{T_k^1}\Big)\setminus\{(f^{2n-1}(\rho_1),
...,\rho_1,0,...)\}\qquad \quad n>0,
$$

\noindent Now, it is enough to notice that
$$
\Theta\Big(O_n(\w_1)\setminus \{\w_1\}\Big)=(a_\infty,-\sum_{k=1}^{n-1}d_{2k-1}).
$$

\noindent Thus puting  $\Theta(\w_1)=a_\infty$ we obtain a  homeomorphism from
$\X_{\omega_0}\cup
\{\w_1\}$ onto $[a_\infty,b_\infty)$.

\hfill$\blacksquare$

In view of  Proposition\,4 the mapping $\Psi$ from Proposition\,12
can be extended to as act bijectively from $\X$ onto  $(a,1]\cup[a_\infty,b_\infty]\times
\{\infty\}$ where
$$\Psi(\x)=(\Theta(\x),\infty), \qquad \x\in \X\setminus\X_1.
$$

\noindent Thus the only thing left for us  to do is to establish  the topology on
$(a,1]\cup[a_\infty,b_\infty]\times \{\infty\}$ inherited from $\X$ by $\Psi$, and
the action of the conjugated mapping $\Psi \circ \f \circ \Psi^{-1}$.

\begin{thm}\label{rho_1>rhotheorem}
Let $(X,f)$ be as in case three a) of Lemma\,1.  Then up to  topological conjugacy
$$\X = (a,1]\cup [a_\infty,b_\infty]\times \{\infty\},$$

\noindent where the interval $(a,1]$ is open in $\X$, and  set
$[a_\infty,b_\infty]\times \{\infty\}$ is naturally homeomorphic
to the  closed interval $[a_\infty,b_\infty]$, see Fig.\,7. A base of
neighborhoods $\{O_{N,\varepsilon}(x,\infty)\}_{N,\varepsilon >
0}$ of a point $(x,\infty)$ in  $[a_\infty,b_\infty]\times
\{\infty\}$ is described as follows:

1) If $x\in (a_\infty,b_\infty)$,  then $O_{N,\varepsilon}(x,\infty)$
consists of the set $\{(y,\infty): y\in (x-\varepsilon,x+\varepsilon)\}$
  and  intervals  $(x_{n,\pm}-\varepsilon, x_{n,\pm}+\varepsilon)\subset (a,1]$ for  $n > N$,
where
\begin{equation}\label{x_plus}
x_{n,+}=x - \Big(2nd_0 + 2\sum_{k=1}^{n-1}(n-k)(d_{2k-1} + d_{2k})\Big),
\end{equation}
\begin{equation}\label{x_minus}
x_{n,-}=- x- \Big(2nd_0 + 2\sum_{k=1}^{n-1}(n-k)(d_{2k-1} + d_{2k}) +
2\sum_{k=1}^{n}d_{2k-1}\Big).
\end{equation}

2) If $x=a_\infty$, then $O_{N,\varepsilon}(a_\infty,\infty)$ contains a set $D_{N,w_1}$
consisting of
points of the form $(-1)^n w_1+2\sum_{k=1}^{n-1}(-1)^k f^{2k-1}(\rho_1)$, for all  $n>N$,
and
$$
O_{N,\varepsilon}(a_\infty,\infty)=\bigcup_{y \in D_{N,w_1}} O_{N,\varepsilon}(y,\infty)\cup
\{(a_\infty,\infty)\}.
$$

3) If  $x= b_\infty$, then $O_{N,\varepsilon}(a_\infty,\infty)$ contains a set $D_{N,w_2}$
consisting of
points of the form $(-1)^n w_2+2\sum_{k=0}^{n-1}(-1)^k f^{2k}(\rho_1)$, for all  $n>N$, and
$$
O_{N,\varepsilon}(b_\infty,\infty)=\bigcup_{y \in D_{N,w_2}} O_{N,\varepsilon}(y,\infty)\cup
\{(b_\infty,\infty)\}.
$$

\noindent Moreover, the induced homeomorphism $\f$ acts as follows:

a) It maps interval $[a_\infty,b_\infty]\times\{\infty\}$ onto
itself in such a manner  that  interval of length $d_0$, marked on
Fig.7, is mapped onto sum of intervals of length $d_0$ and $d_1$;
interval of length $d_{2n}$, $n>0$,  is mapped onto  interval of
length $d_{2n+1}$, and  interval of length $d_{2n-1}$  is mapped
onto  interval of length $d_{2n}$. There is a  fixed point
$(\omega_0,\infty)$ and a cycle $\{(a,\infty), (b,\infty)\}$ of
period two.

b) It is an increasing bijection of $(a,1]$ onto itself with fixed point $1$:  $\f$ maps
$(0,1]$ onto $(-d_0,1]$,
and  then (looking from right -- see Fig.\,7) $\f$ maps every ``climb'' onto the next ``descent''
and  every ``descent''
onto the next ``climb''. It maps interval of length $d_{2n}$, $n>0$,   onto  interval of length
$d_{2n+1}$, and
 interval of length $d_{2n-1}$   onto  interval of length $d_{2n}$.
\end{thm}

\begin{center}
\setlength{\unitlength}{0.9mm}
\begin{picture}(120,90)(-100,-41)
\thinlines
\qbezier(0,3)(0,3)(25,20)\put(25,20){\circle*{1}}\put(25,21){\scriptsize{$1$}}
\qbezier(0,3)(-6,10.5)(-12,18)\put(-0,3){\circle*{1}}\put(0,0.5){\scriptsize{$0$}}
\qbezier(-24,2)(-18,10)(-12,18)\put(-12,18){\circle*{1}}\put(-13,19){\scriptsize{-$d_0$}}
\qbezier(-24,2)(-24,2)(-28,-7)\put(-24,2){\circle*{1}}\put(-23,1){\scriptsize{-$2d_0$}}
\qbezier(-32,2)(-32,2)(-28,-7)\put(-28,-7){\circle*{1}}\put(-32,-9.5){\scriptsize{-$2d_0$-$d_1$}}
\qbezier(-32,2)(-32,2)(-42,18.5)\put(-32,2){\circle*{1}}\put(-41.5,2){\scriptsize{-$2d_0$-$2d_1$}}
\qbezier(-45,28)(-45,28)(-42,18.5)\put(-42,18.5){\circle*{1}}\put(-41,18){\scriptsize{-$3d_0$-$2d_1$}}
\qbezier(-45,28)(-45,28)(-48,18.5)\put(-45,28){\circle*{1}}\put(-52,29){\scriptsize{-$3d_0$-$2d_1$-$d_2$}}
\qbezier(-56,1)(-56,1)(-48,18.5)\put(-48,18.5){\circle*{1}}\put(-63,18){\scriptsize{-$3d_0$-$2d_1$-$2d_2$}}
\qbezier(-59,-8.5)(-59,-8.5)(-56,1)\put(-56,1){\circle*{1}}\put(-55.5,-0.5){\scriptsize{-$4d_0$-$2d_1$-$2d_2$}}
\qbezier(-59,-8.5)(-59,-8.5)(-60.5,-14)\put(-59,-8.5){\circle*{1}}\put(-57.5,-9){\scriptsize{-$4d_0$-$3d_1$-$2d_2$}}
\qbezier(-62,-8.5)(-62,-8.5)(-60.5,-14)\put(-60.5,-14){\circle*{1}}\put(-69,-17){\scriptsize{-$4d_0$-$3d_1$-$2d_2$-$d_3$}}
\qbezier(-62,-8.5)(-62,-8.5)(-64,1.5)\put(-62,-8.5){\circle*{1}}\put(-76,-7){\tiny{-$4d_0$-$3d_1$-$2d_2$-$2d_3$}}
\qbezier(-70,19.5)(-70,19.5)(-64,1.5)\put(-64,1.5){\circle*{1}}\put(-76,3){\tiny{-$4d_0$-$4d_1$-$2d_2$-$2d_3$}}
\qbezier(-70,19.5)(-70,19.5)(-72,29)\put(-70,19.5){\circle*{1}}
\qbezier(-73.5,35)(-73.5,35)(-72,29)\put(-72,29){\circle*{1}}

%\qbezier(-74.5,38.5)(-74.5,38.5)(-73.5,35)\put(-73.5,35){\circle*{1}}
%\qbezier(-74.5,38.5)(-74.5,38.5)(-75.5,35)\put(-75.5,35){\circle*{1}}
\qbezier(-74.5,29)(-74.5,29)(-73.5,35)\put(-73.6,35){\circle*{1}}
\qbezier(-74.5,29)(-74.5,29)(-76,19)\put(-74.5,29){\circle*{1}}
\qbezier(-80,0.5)(-80,0.5)(-76,19.5)\put(-76,19.5){\circle*{1}}
\qbezier(-80,0.5)(-80,0.5)(-81.5,-9.5)\put(-80,0.5){\circle*{1}}
\qbezier(-82.5,-15.5)(-82.5,-15.5)(-81.5,-9.5)\put(-81.5,-9.5){\circle*{1}}
\qbezier(-82.5,-15.5)(-82.5,-15.5)(-83.25,-19)\put(-82.5,-15.5){\circle*{1}}
\qbezier(-84,-15.5)(-84,-15.5)(-83.25,-19)\put(-83.25,-19){\circle*{1}}
\qbezier(-84,-15.5)(-84,-15.5)(-85,-9.5)
\qbezier(-86,0.5)(-86,0.5)(-85,-9.5)
\qbezier(-86,0.5)(-86,0.5)(-88,20)
\qbezier(-89,30)(-89,30)(-88,20)
\qbezier(-89,30)(-89,30)(-89.75,36)
\qbezier(-90.25,40)(-90.25,40)(-89.75,36)
\qbezier(-90.25,40)(-90.25,40)(-91,36)
\qbezier(-94,-22)(-92.5,8)(-91,36)
\qbezier(-94,-22)(-95.75,10)(-96.5,42)
\qbezier[220](-98.5,-23)(-97.5,9.5)(-96.5,42)
\qbezier[200](-98.5,-23)(-99.25,10)(-100,43)
\qbezier[190](-101,-24)(-100.5,9.5)(-100,43)
\qbezier[170](-101,-24)(-101,10)(-102,44)
\qbezier[140](-102,-24.5)(-102,10)(-102,44)

\qbezier[100](10,-20)(-50,-20)(-103,-20)
\qbezier[150](10,-16)(-50,-16)(-103,-16)
\qbezier[200](10,-10)(-50,-10)(-103,-10)
\qbezier[300](10,0)(-50,0)(-103,0)
\qbezier[300](10,20)(-50,20)(-103,20)
\qbezier[200](10,30)(-50,30)(-103,30)
\qbezier[150](10,36)(-50,36)(-103,36)
\qbezier[150](10,40)(-50,40)(-103,40)

\qbezier(-103,-25)(-103,10)(-103,45)
\put(-103,-25){\circle*{1}}\put(-109,-28){\footnotesize{$(a_\infty,\infty)$}}
\put(-103,45){\circle*{1}}\put(-109,47){\footnotesize{$(b_\infty,\infty)$}}

\put(-106.5,9){\footnotesize{$d_0$}}

\put(-106.5,-6){\footnotesize{$d_1$}}
\put(-106.5,24){\footnotesize{$d_2$}}
\put(-106.5,-14){\footnotesize{$d_3$}}
\put(-106.5,32){\footnotesize{$d_4$}}
\put(-106.5,-19){\footnotesize{$d_5$}}
\put(-106.5,37){\footnotesize{$d_6$}}
\put(-106,-24){\footnotesize{$\vdots$}}
\put(-106,40){\footnotesize{$\vdots$}}

\put(-45,-30){\small{Fig. 7.}}
\put(-100,-37){\small{The inverse limit associated with $([0,1],f)$ in case $f$ is
$\F_{2^1}$ mapping}}
\put(-100,-42){\small{possessing an
attractive cycle of period two.}}
\end{picture}

\end{center}

\vspace{3mm}

\noindent {\it Proof.}
In view of Propositions\,4 and \ref{sickhomeomorphism},
to describe the topology on $\X$ we only need to consider  neighborhoods of a point
$(x,\infty)\in [a_\infty,b_\infty]\times \{\infty\}$.
We have the following cases:
\Item{1)} $x \in (a_\infty,b_\infty)$. We assume that $x<0$, the case $x\geq 0$ is similar.
Then   $x$ belongs to $[-\sum_{k=1}^{N_x} d_{2k-1},-\sum_{k=1}^{N_{x}-1} d_{2k-1})$, for
some $N_x\in \N$,
and there exists $\x =(x_0,x_1,...)\in \X_{T^1_{N_x}}$ such that $\Theta(\x)=x$ ($x_0=x$).

Consider the base of neighborhoods $\{O_{N,\epsilon}(\x)\}_{N > 2N_x,\epsilon >0}$ of
the point $\x=(x_0,x_1,...)$ where
$$
O_{N,\epsilon}(\x)=\{(y_0,y_1,...)\in\X: y_{2N_x-1}\in(\rho,1),\,\,\, |y_{N} -
x_{N}|<\epsilon\}.
$$

\noindent For each pair $(N,\epsilon)$ there exists $\varepsilon> 0$ such that
$(x_{0}-\varepsilon,x_{0} +\varepsilon)\subset f^{N}(x_{N}-\epsilon,x_{N} +\epsilon)$,
and then
 $V_{N,\varepsilon}(\x)=\{(y_0,y_1,...)\in O_{N,\epsilon}(\x): |x_0-y_0|<\varepsilon\}$
 is open neighborhood of $\x$ and  $V_{N,\varepsilon}(\x)$ (for sufficiently small
 $\epsilon$) has the following  form
$$
V_{N,\varepsilon}(\x)=\{(y_0,y_1,...)\in \bigcup_{n=N}^{\infty}\,\X_{T_{n,N}^{(1)}}\cup
\X_{T_{N_x}^{1}}: y_0\in
 (x_{0}-\varepsilon,x_{0} +\varepsilon)\}\subset O_{N,\epsilon}(\x).
$$

Mapping $\Psi$ carries points $\y_{n,-}=(x_0,...)\in \X_{T_{2n,N}^{(1)}}$ and $\y_{n,+}=
(x_0,...)\in \X_{T_{2n+1,N}^{(1)}}$
 into the points $x_{n,-}$ and $x_{n,+}$ given by (\ref{x_plus}),(\ref{x_minus}), cf.
 Proposition\,3, and
 $$
 \Psi(\{\y\in\X_{T_{N_x}^{1}}: y_0\in  (x_{0}-\varepsilon,x_{0} +\varepsilon)\})=
 \{(y,\infty): y\in (x-\varepsilon,x+\varepsilon)\}.
 $$

\noindent Hence $\Psi(V_{N,\varepsilon}(\x))=O_{N,\varepsilon}(x,\infty)$ and  the
hypotheses of the
  item 1) follows.
\Item{2)} $x = a_\infty$. Then $\Psi(\w_1)=(x,\infty)$. Consider the base of neighborhoods
$\{O_{N,\epsilon}(\w_1)\}_{N ,\epsilon >0}$ of   $\w_1=(w_1,w_2,w_1,...)$ where
$$
O_{N,\epsilon}(\w_1)=\{(y_0,y_1,...)\in\X: y_{2N}\in (w_1-\epsilon,w_1 +\epsilon)\}.
$$

\noindent For each pair $(N,\epsilon)$ there exists $\varepsilon> 0$ such that
$(w_1-\varepsilon,w_1 +\varepsilon)\subset f^{2N}(w_1-\epsilon,w_1 +\epsilon)$, and then
$V_{N,\varepsilon}(\w_1)=\{(y_0,y_1,...)\in O_{N,\epsilon}(\w_1): |\omega_1-y_0|
<\varepsilon\}$ is open
neighborhood of $\w_1$ and $V_{N,\varepsilon}(\w_1)$ (for sufficiently small $\epsilon$)
has the following  form
$$
V_{N,\varepsilon}(\w_1)=\{(y_0,y_1,...)\in \bigcup_{n=2N}^{\infty}\,\X_{T_{n,N}^{(1)}}\cup
\X_{T_{n}^{1}}\cup \X_{T_\infty^1}: y_0\in  (w_1-\varepsilon,w_1 +\varepsilon)\}.
$$

\noindent In an analogous fashion as in the proof of the item 1) we conclude
that the image of $V_{N,\varepsilon}(\w_1)$ under $\Psi$ is
$O_{N,\varepsilon}(a_\infty,\infty)$ and  the hypotheses of  the
item 2) follows. \Item{3)} $x = b_\infty$. The proof is analogous
to the proof of the item 2) and hence omitted.

The item a) follows immediately from Propositions
\ref{X_omega_in_case_3}, \ref{X minus X_1 in case_3a)} and fact
that $\f(\w_1)=\w_2$, $\f(\w_2)=\w_1$.

To show up item b) we notice that $\f(\X_{T_0})=\X_{T_0}\cup \X_{T_{1}}$ and
$$
\f(\X_{T_n})=\X_{T_n}\cup \X_{T_{n+1,1}^{(1)}},\qquad \f(\X_{T_{n,k}^{(1)}})=
\X_{T_{n+1,k}^{(0)}}
 ,\qquad \f(\X_{T_{n,k}^{(0)}})=\X_{T_{n+1,k+1}^{(1)}},
$$

\noindent for $n>0$, and $k=1,...,\Big[\frac{n-1+i}{2}\Big]$. In view of the form of
the homeomorphism
constructed in subsection 3.1 the hypotheses now easily follows.

\hfill$\blacksquare$

\subsubsection{The case of two intervals of periodic points }

In the case 3b) from the Lemma\,1 the points from sets $[0,f(\rho_1)]
\setminus\{w_1\}$ and
$[\rho,\rho_1]\setminus\{w_2\}$  are periodic points of period four and set $\{w_1,w_2\}$
forms a cycle of period two. Then we have
$$
d_1= d_3 = d_5 = ...=f(\rho_1)\, ,\quad\,  \qquad d_2 = d_4 = d_6 = ...=\rho_1-\rho\,\qquad
\text{and}
$$
$$\qquad \Phi(\X_{T_{\infty}^1})=[0,f(\rho_1)], \quad  \Phi(\X_{T_{\infty}^0})=[\rho,\rho_1].
$$

\noindent Thus in Proposition\,3 we have  $a_\infty=-\infty$ and
$b_\infty=+\infty$.
 We prolong now the homeomorphism from the Proposition\,3 to as act
 from the set
   $\X\setminus \X_1=\X_{\omega_0}\cup\X_{T_{\infty}^0}\cup\X_{T_{\infty}^1}$ onto a certain
   topological space,
    cf. Proposition\,4.

\begin{prop}\label{X minus X_1 in case_3b)}
Let $(X,f)$ be as in case 3b) in Lemma\,1, and let $\Theta$ be the mapping from
Proposition\,3.
If we define
$$
\Theta(\x)=(\Phi(\x),-\infty) ,\textrm{ for } \x\in\X_{T_{\infty}^1}  \,\,\,\,
\textrm{ and }\,\,\,\,
 \Theta(\x)=(\Phi(\x),+\infty),\textrm{ for } \x\in\X_{T_{\infty}^0},
$$

\noindent then $\Theta$ is a homeomorphism from $\X\setminus \X_1$ onto a topological space
 $$[0,f(\rho_1)]\times\{-\infty\}\cup(-\infty,+\infty)\cup  [\rho,\rho_1]\times\{+\infty\} $$

\noindent consisting of the real axis and two closed subsets homeomorphic to closed intervals
(see Fig.\,8).
  The base of neighborhoods  $\{O_{N,\varepsilon}(x,\pm\infty)\}_{N,\varepsilon > 0}$ of a
  point $(x,\pm\infty)$
  is described as follows:

\item{1)} $O_{N,\varepsilon}(x,-\infty)$  consists of the set
$\{(y,-\infty): y\in (x-\varepsilon,x+\varepsilon)\cap [0,d_1]\}$ and
intervals, for all $n > N$,
\begin{equation}\label{neighbor1}
((-1)^{n} x -
2\Big[\frac{n}{2}\Big]d_1 -\varepsilon,(-1)^{n} x -
2\Big[\frac{n}{2}\Big]d_1 +\varepsilon) \subset (-\infty,0).
\end{equation}

\item{2)} $O_{N,\varepsilon}(x,-\infty)$  consists of the set
$\{(y,+\infty): y\in (x-\varepsilon,x+\varepsilon)\cap [0,d_2]\}$ and
intervals, for all $n > N$,
$$
((-1)^{n+1} (x-\rho) +
2\Big[\frac{n}{2}\Big]d_2 +\rho-\varepsilon,(-1)^{n+1} (x-\rho) +
2\Big[\frac{n}{2}\Big]d_2 + \rho +\varepsilon) \subset (\rho,+\infty). $$

\noindent The conjugated mapping $\Theta\circ  \f \circ \Theta^{-1}$ is a
decreasing  homeomorphism of  $(a_\infty,b_\infty)$ such that for
all $n>0$,  $ (-n d_1,-(n-1)d_1)$ is mapped onto $(\rho+n d_2,
\rho + (n+1)d_2)$  and  $(\rho+n d_2,\rho+  (n+1)d_2)$ is mapped
onto $(-(n+1) d_1,-n d_1)$, cf. Fig. 8. Moreover we have
$$
\Theta\circ  \f \circ \Theta^{-1}(x,\pm \infty)=(f(x),\mp\infty),
$$

\noindent that is  $[0,f(\rho_1)]\times\{-\infty\}$ and $[\rho,\rho_1]\times\{+\infty\} $
are intervals of periodic points with period four.
\end{prop}

\begin{center}
\setlength{\unitlength}{0.6mm}
\begin{picture}(220,59)(-90,-19)
\thinlines

\qbezier(35,20)(35,20)(63,32)\put(35,20){\circle*{1.25}}\put(33,21.5){\scriptsize{$\rho$}}
\put(57,34){\scriptsize{$\rho_1=\rho+d_2$}}
\qbezier(63,32)(63,32)(80,6.5)\put(63,32){\circle*{1.25}}\put(75,3){\scriptsize{$\rho$+$2d_2$}}
\qbezier(80,6.5)(80,6.5)(92,34)\put(80,6.5){\circle*{1.25}}\put(86.5,36){\scriptsize{$\rho$+$3d_2$}}
\qbezier(92,34)(92,34)(99,5)\put(92,34){\circle*{1.25}}
\qbezier(99,5)(99,5)(103,35)
\qbezier(103,35)(103,35)(105.5,5)
\qbezier(105.5,5)(105.5,5)(107,35)
\qbezier[150](107,35)(107,35)(108,5)
\qbezier[100](108,5)(108,5)(109,35)
\qbezier[75](109,35)(109,20)(109,5)

\qbezier(110,5)(110,20)(110,35)
\put(110,5){\circle*{1.25}}\put(104,1){\scriptsize{$(\rho,$+$\infty)$}}
\put(110,35){\circle*{1.25}}\put(104,37){\scriptsize{$(\rho_1,$+$\infty)$}}

%\qbezier(-4,7.5)(9.5,18)(23,28.5)\put(23,28.5){\circle*{1}}\put(24,29){$1$}
\qbezier(-5,20)(15,20)(35,20)\put(-5,20){\circle*{1.25}}\put(-7,21.5){\scriptsize{$0$}}
%\qbezier(-5,20)(-5,20)(-32,7)
\qbezier(-5,20)(-33,32)(-33,32)
%\put(-32,7){\circle*{1.25}}\put(-31,4.5){-$2\rho$}
\qbezier(-33,32)(-49,6.5)(-49,6.5)\put(-33,32){\circle*{1.25}}\put(-36,34){\scriptsize{-$d_1$}}
\qbezier(-49,6.5)(-49,6.5)(-61,34)\put(-49,6.5){\circle*{1.25}}\put(-53,3){\scriptsize{-$2d_1$}}
\qbezier(-61,34)(-61,34)(-69,5)\put(-61,34){\circle*{1.25}}\put(-64,36){\scriptsize{-$3d_1$}}
\qbezier(-69,5)(-71,20)(-73,35)
\qbezier(-75.5,5)(-73,35)(-73,35)
\qbezier(-75.5,5)(-75.5,5)(-77,35)
\qbezier[150](-78,5)(-77.5,20)(-77,35)
\qbezier[100](-78,5)(-78.5,20)(-79,35)
\qbezier[75](-79,5)(-79,20)(-79,35)

\qbezier(-80,5)(-80,20)(-80,35)
\put(-80,5){\circle*{1.25}}\put(-86,1){\scriptsize{$(0,-\infty)$}}
\put(-80,35){\circle*{1.25}}\put(-88,37){\scriptsize{$(f(\rho_1),-\infty)$}}

\put(0,-5){\small{Fig. 8.}}
\put(-65,-14){\small{The subspace $\X\setminus \X_1$ of the inverse limit associated with $([0,1],f)$ }}
\put(-65,-21){\small{in case $f$ is $\F_{2^2}$ mapping.}}
\end{picture}
\end{center}
\vspace{3mm}

\noindent {\it Proof.}
We shall only consider  the neighborhoods from the item  1), the part of proof concerning
the item 2) is similar.\\
Let $\x =(x_0,x_1,...)\in \X_{T_{\infty}^1}$ be such that $\Phi(\x)=x$, that is $x_0=x$.
The sets
$O_{N,\varepsilon}(\x)=\{(y_0,...)\in \X\setminus \X_1: y_{4n}\in(x_{4n}-\varepsilon,x_{4n}
+\varepsilon)\cap (0,f(\rho_1))\}$
form  a base of neighborhoods of $\x$ in $\X\setminus \X_1$. As $(0,f(\rho_1))$ consists of
periodic points of period four  we have
$$
O_{N,\varepsilon}(\x)=\{(y_0,y_1,...)\in \bigcup_{k=N}^\infty \X_{T_k^1}\cup
\X_{T_{\infty}^{1}}: y_0\in
 (x_{0}-\varepsilon,x_{0} +\varepsilon)\}.
$$

\noindent It is clear that
$$\Theta(\{(y_0,...)\in \X_{T_{\infty}^{1}}: |x_{0}-y_{0}|< \varepsilon\})=\{(y,-\infty):
|x-y|<\varepsilon,\, y\in [0,d_1]\}.$$

\noindent  For each point $\y=(y_0,y_1,...)\in \X_{T_n^1}$ such that $y_0=x_0=x$ we have
by formula (\ref{thirdformula}) that
$\Theta(\y)=(-1)^n x+2\sum_{k=1}^{n-1}(-1)^k f^{2k-1}(\rho_1)=(-1)^{n} x - 2[\frac{n}{2}]d_1$.
Thus we have $\Theta(O_{N,\varepsilon}(\x))=O_{N,\varepsilon}(x,-\infty)$, that is  $\Theta$
is a
 homeomorphism from $\X_{\omega_0}\cup \X_{T_{\infty}^1}$ onto $[0,f(\rho_1)]
 \times\{-\infty\}\cup(-\infty,+\infty)$.
By the same argumentation concerning now elements from
$\X_{T_{\infty}^0}$ one gets  $\Theta$ as a homeomorphism
acting  on $\X\setminus \X_1$.
\\
The form of the conjugated mapping follows from Proposition \ref{X_omega_in_case_3} and
from the fact that  $\f$ maps interval  $[0,f(\rho_1)]$ onto $[\rho, \rho_1]$.

\hfill$\blacksquare$

We are now in a position to describe homeomorphic image of $\X$ in the case 3b) from the
Lemma\,1.
 By virtue of the Proposition\,3 the mapping $\Psi$ from the
 Proposition
 \ref{sickhomeomorphism} can be extended to as act bijectively from $\X$ onto  $[0,f(\rho_1)]
 \times\{-\infty\}
 \cup\R\times \{\infty\}\cup  [\rho,\rho_1]\times\{+\infty\}\cup (-\infty,1]$ where
$$\Psi(\x)=(\Theta(\x),\infty), \qquad \x\in \X_{\omega_0}, \qquad \qquad \Psi(\x)=\Theta(\x),
\qquad \x\in \X_{T_{\infty}^{i}}.
$$
Let us  describe now in detaile the topology on $\Big([0,f(\rho_1)]\times\{-\infty\}\Big)\cup
\Big(\R\times \{\infty\}\Big)\cup
\Big([\rho,\rho_1]\times\{+\infty\}\Big)\cup(-\infty,1]$ inherited from $\X$.

\begin{thm}\label{rho_1>rhotheoremb)}
Let $(X,f)$ be as in the case three b) of the Lemma\,1.  Then up to  topological
conjugacy
$$\X=\underbrace{\Big([0,f(\rho_1)]\times\{-\infty\}\Big)}_{\X_{T_\infty^1}}
\cup\underbrace{\Big(\R\times \{\infty\}\Big)}_{\X_{\omega_0}}\cup
 \underbrace{\Big([\rho,\rho_1]\times\{+\infty\}\Big)}_{\X_{T_\infty^2}}\cup
 \underbrace{(-\infty,1]}_{\X_{1}},$$

\noindent where the topology on $\X$ is described as follows (see Fig.\,9):

\Item{0)} The interval $(-\infty,1]$ is open in $\X$.

\Item{1)} A point $(x,\infty)$ in $\R\times \{\infty\}$ has a base of neighborhoods
$\{O_{N,\varepsilon}(x,\infty)\}_{N,\varepsilon > 0}$ where
$O_{N,\varepsilon}(x,\infty)$ consists of the set $\{(y,\infty): y\in (x-\varepsilon,x+
\varepsilon)\}$
 and  all intervals  $(x_{n,\pm}-\varepsilon, x_{n,\pm}+\varepsilon)$ for  $n > N$,
where
\begin{equation}\label{x_plus'}
x_{n,+}=x - \Big(2nd_0 + n(n-1)(d_{1} + d_{2})\Big),
\end{equation}
\begin{equation}\label{x_minus'}
x_{n,-}=- x- \Big(2nd_0 + n^2d_{1} + n(n-1)d_{2})\Big).
\end{equation}

\Item{2)} A point $(x,-\infty)$ in $[0,f(\rho_1)]\times\{-\infty\}$ has a base of
neighborhoods
$\{O_{N,\varepsilon}(x,-\infty)\}_{N,\varepsilon > 0}$ where $
O_{N,\varepsilon}(x,-\infty)$ consists of the set $\{(y,-\infty): y\in (x-\varepsilon,x
+\varepsilon)\cap [0,d_1]\}$ and  the sum
$$
\bigcup_{n=N}^{\infty}O_{n,\varepsilon}\Big((-1)^{n} x - 2\Big[\frac{n}{2}\Big]d_1,\infty\Big).
$$
\Item{3)} A point $(x,+\infty)$ in $[\rho,\rho_1]\times\{+\infty\}$ possess a base of
neighborhoods
$\{O_{N,\varepsilon}(x,+\infty)\}_{N,\varepsilon > 0}$ where $
O_{N,\varepsilon}(x,+\infty)$ consists of the set $\{(y,+\infty): y\in (x-\varepsilon,x
+\varepsilon)\cap [\rho,\rho_1]\}$
and  the sum
$$
\bigcup_{n=N}^{\infty}O_{n,\varepsilon}\Big((-1)^{n+1} (x-\rho) + 2\Big[\frac{n}{2}\Big]d_2
+\rho,\infty\Big).
$$
Moreover, the induced homeomorphism $\f$ acts as follows:
\par
a) It is a decreasing bijection of $(a,1]$ onto itself with fixed point $1$, see Fig.\,8, cf.
Theorem\,4.
\par
b) It is a ``decreasing'' bijection of $(a_\infty,b_\infty)\times\{\infty\}$ onto itself with
fixed point
$(\omega_0,\infty)$, cf. Proposition\,5.
\par
 c) It maps $[0,f(\rho_1)]\times\{-\infty\}$ onto $[\rho,\rho_1]\times\{+\infty\} $ and vice versa. Furthermore we have
$$
\f (x,\pm \infty)=(f(x),\mp\infty),
$$

\noindent that is  $[0,f(\rho_1)]\times\{-\infty\}$ and $[\rho,\rho_1]\times\{+\infty\} $
are intervals of
periodic points with period four, and $\{(w_1,-\infty),(w_2,+\infty)\}$ is a cycle of period
two.
\end{thm}

\vspace{3mm}

\noindent {\it Proof.}
 In view of  Propositions\,5 and \ref{sickhomeomorphism}, we only
 need to consider
 neighborhoods of  points: $(x,\infty)\in \R\times \{\infty\}$,   $(x,-\infty)\in
 [0,f(\rho_1)]\times\{-\infty\}$
 and $(x,+\infty)\in [\rho,\rho_1]\times\{+\infty\}$:

\Item{1)} Let $(x,\infty)\in \R\times \{\infty\}$. We arrive at the
hypotheses by the same argumentation as in the  proof of  the item
1) from Theorem\,4. Recall that
 $d_1= d_3 = d_5 = ...$, and  $d_2 = d_4 = d_6 = ...$,
  hence formulas (\ref{x_plus}) and  (\ref{x_minus}) reduce to  (\ref{x_plus'}) and
  (\ref{x_minus'}).

\Item{2)} If $(x,-\infty) \in [0,f(\rho_1)]\times\{-\infty\}$, then there exists
$\x\in \X_{T_\infty^1}$
such that $\Psi(x)=(x,\infty)$. Hence $\x=(x_0,x_1,x_2,...)$ where $x_n=x_{n+4}$
for all $n\in\N$.

Consider a base of neighborhoods $\{O_{N,\epsilon}(\x)\}_{N > 0}$  where
$$
O_{N,\epsilon}(\x)=\{(y_0,y_1,...)\in\X: y_{4N}\in (x_0-\epsilon,x_0 +\epsilon)\}.
$$

\noindent For each pair $(N,\epsilon)$ there exists $\varepsilon> 0$ such that
 $(x_0-\varepsilon,x_0 +\varepsilon)\cap [0,f(\rho_1)]$ is contained in $f^{4N}(x_0
 -\epsilon,x_0 +\epsilon)$.
 Then  $O_{N,\epsilon}(\x)$ contains the open set $V_{N,\varepsilon}(\x)=\{(y_0,y_1,...)
 \in O_{N,\epsilon}(\x):
  y_0\in  (x_0-\varepsilon,x_0 +\varepsilon)\cap [0,f(\rho_1)]\}$ and we have
$$
V_{N,\varepsilon}(\x)=\{(y_0,...)\in \bigcup_{n=4N}^{\infty}\X_{T_{n,N}^{(1)}}\cup
\X_{T_{n}^{1}}\cup \X_{T_\infty^1}:
 y_0\in  (x_0-\varepsilon,x_0 +\varepsilon)\cap [0,f(\rho_1)]\}.
$$

\noindent In view of the form of homeomorphisms from Propositions \ref{sickhomeomorphism},
5
one can see that $\Psi(V_{N,\varepsilon}(\x))=O_{N,\varepsilon}(x,-\infty)$ and  the
hypotheses of  the item 2) follows.
\Item{3)} Let $(x,+\infty)\in [\rho,\rho_1]\times\{+\infty\}$. The proof is similar to
the proof of the item\,\,2)
and hence it is omitted.

The form of the induced homeomorphism $\f$  follows from Propositions\,2
and  5.

\hfill$\blacksquare$

\subsection{The $C^*$-algebraic context}\label{C^*}\label{6}

In \cite{maxid} the author and his co-author came across the
construction discussed in this paper  while investigating the
spectrum of a certain commutative $C^*$-algebra, see also
\cite{covalg}. For the sake of completness  we include  a brief
discussion concerning  this context.

Let us consider a unitary operator $U$ and a unital commutative $C^*$-algebra $\A$
acting nondegenerately
 on some Hilbert space $\H$.
 Moreover, let us assume that
 $$U \A U^*\subset \A \,\, \textrm{ but }\,\, U^* \A U\nsubseteq \A.$$

\begin{center}
\setlength{\unitlength}{0.9mm}
\begin{picture}(100,220)(10,-30)

\thinlines

% Sprężynka

\put(63.7,169.7){\circle*{1}}\put(61,169.2){\footnotesize{$1$}}
\qbezier(63.6,169.7)(63.6,169.7)(55.8,132.5)
\put(55.8,132.5){\circle*{1}}\put(53.1,132.2){\footnotesize{$0$}}
\qbezier(68.8,128)(68.8,128)(55.8,132.5)
\put(68.8,128){\circle*{1}}\put(68.4,129){\footnotesize{-$d_0$}}
\qbezier(68.8,128)(68.8,128)(54.9,108)
\put(54.9,108){\circle*{1}}\put(50,109.5){\footnotesize{-$2d_0$}}
\qbezier(39,100.3)(39,100.3)(54.9,108)
\put(39,100.3){\circle*{1}}\put(32,102){\footnotesize{-$2d_0$-$d_1$}}
\qbezier(39,100.3)(39,100.3)(54.9,100)
\put(54.9,100){\circle*{1}}\put(51,101.6){\scriptsize{-$2d_0$-$2d_1$}}
\qbezier(70.05,98.6)(70.05,98.6)(54.9,100)
\put(70.05,98.6){\circle*{1}}\put(65,99.8){\scriptsize{-$3d_0$-$2d_1$}}
\qbezier(70.05,98.6)(70.05,98.6)(75.9,105)
\put(75.9,105){\circle*{1}}\put(70,106){\scriptsize{-$3d_0$-$2d_1$-$d_2$}}
\qbezier(70,92.6)(70,92.6)(75.9,105)
\put(70,92.6){\circle*{1}}\put(63,90){\scriptsize{-$3d_0$-$2d_1$-$2d_2$}}
\qbezier(70,92.6)(70,92.6)(54,75.9)
\put(54,75.9){\circle*{1}}\put(48,77){\scriptsize{-$4d_0$-$2d_1$-$2d_2$}}
\qbezier(37.1,69.6)(37.1,69.6)(54,75.9)
\put(37.1,69.6){\circle*{1}}
\qbezier(37.1,69.6)(37.1,69.6)(40.8,57.7)
\put(40.8,57.7){\circle*{1}}
\qbezier(37.1,66.54)(37.1,66.54)(40.8,57.7)
\put(37.1,66.54){\circle*{1}}
\qbezier(37.1,66.54)(37.1,66.54)(53.85,67.7)
\put(53.85,67.7){\circle*{1}}
\qbezier(70.2,71.2)(70.2,71.2)(53.85,67.7)
\put(70.2,71.2){\circle*{1}}
\qbezier(70.2,71.2)(70.2,71.2)(76.3,78.2)
\put(76.3,78.2){\circle*{1}}
\qbezier(94.55,74.7)(94.55,74.7)(76.3,78.2)
\put(94.55,74.7){\circle*{1}}
\qbezier(94.55,74.7)(94.55,74.7)(76.35,75.28)
\put(76.35,75.28){\circle*{1}}
\qbezier(70.25,63.2)(70.25,63.2)(76.35,75.28)
\put(70.25,63.2){\circle*{1}}
\qbezier(70.25,63.2)(70.25,63.2)(53.5,49.5)
\put(53.5,49.5){\circle*{1}}
\qbezier(36.5,43.75)(36.5,43.75)(53.5,49.5)
\put(36.5,43.75){\circle*{1}}
\qbezier(36.5,43.75)(36.5,43.75)(40.6,32.28)
\put(40.6,32.28){\circle*{1}}
\qbezier(23.48,35.7)(23.48,35.7)(40.6,32.28)
\put(23.48,35.7){\circle*{0.95}}
\qbezier(23.48,35.7)(23.48,35.7)(40.3,31.45)
\put(40.3,31.45){\circle*{0.9}}
\qbezier(36.5,40.95)(36.5,40.95)(40.3,31.45)
\put(36.5,40.95){\circle*{0.875}}
\qbezier(36.5,40.95)(36.5,40.95)(53.35,43.43)
\put(53.35,43.43){\circle*{0.85}}
\qbezier(70.74,51.48)(70.74,51.48)(53.35,43.43)
\put(70.74,51.48){\circle*{0.825}}
\qbezier(70.74,51.48)(70.74,51.48)(76.8,59.45)
\put(76.8,59.45){\circle*{0.8}}
\qbezier(95,56.5)(95,56.5)(76.8,59.45)
\put(95,56.5){\circle*{0.775}}
\thinlines
\qbezier(95,56.5)(95,56.5)(88.5,65.4)
\put(88.5,65.4){\circle*{0.67}}
\qbezier(95.05,55.1)(95.05,55.1)(88.5,65.4)
\put(95.05,55.1){\circle*{0.725}}
\qbezier(95.05,55.1)(95.05,55.1)(76.9,56.65)
\put(76.9,56.65){\circle*{0.7}}
\qbezier(70.74,46.78)(70.74,46.78)(76.9,56.65)
\put(70.74,46.78){\circle*{0.675}}
\qbezier[110](70.74,46.78)(62.045,41.355)(53.35,35.93)
\put(53.35,35.9){\circle*{0.65}}
\qbezier[90](36,31.5)(44.675,33.7)(53.35,35.9)
\put(36,31.5){\circle*{0.625}}
\qbezier[90](36,31.5)(38.1,26.25)(40.2,21)
\put(40.2,21){\circle*{0.6}}
\qbezier[90](40.2,21)(31.65,22.7)(23.1,24.4)
\put(23.1,24.4){\circle*{0.575}}
\qbezier[85](32.4,15.4)(27.75,19.9)(23.1,24.4)
\put(32.4,15.4){\circle*{0.575}}
\qbezier[81](32.4,15.4)(27.75,19.7)(23.2,24)
\put(23.1,24){\circle*{0.53}}
\qbezier[80](40.2,20)(31.65,22)(23.1,24)
\put(40.2,20){\circle*{0.5}}
\qbezier[78](36,30.2)(38.1,25.1)(40.2,20)
\put(36,30.2){\circle*{0.5}}
\qbezier[75](36,30.2)(44.675,32)(53.35,33.8)
\put(53.35,33.8){\circle*{0.5}}
\qbezier[90](70.74,43.6)(62.045,38.7)(53.35,33.8)
\put(70.74,43.6){\circle*{0.475}}
\qbezier[90](70.74,43.6)(70.74,43.6)(77.2,53)
\put(77.2,53){\circle*{0.45}}
\qbezier[70](95.6,50.8)(86.4,51.9)(77.2,53)
\put(95.6,50.8){\circle*{0.45}}
\qbezier[69](95.6,50.8)(92.3,55.6)(89,60.4)
\put(89,60.4){\circle*{0.4}}
\qbezier[68](104.6,54.6)(96.8,57.5)(89,60.4)
\put(104.6,54.6){\circle*{0.4}}
\qbezier[67](104.6,54.6)(96.8,57.3)(89,60)
\put(89,60){\circle*{0.4}}
\qbezier[69](95.6,49.8)(92.2,54.9)(88.8,60)
\put(95.6,49.8){\circle*{0.4}}
\qbezier[70](95.6,49.8)(86.4,50.8)(77.2,51.8)
\put(77.2,51.8){\circle*{0.4}}
\qbezier[68](70.74,41.8)(73.97,46.8)(77.2,51.8)
\put(70.74,41.8){\circle*{0.4}}
\qbezier[80](70.74,41.8)(62.045,36.7)(53.35,31.6)
\put(53.35,31.6){\circle*{0.4}}
\qbezier[63](36,27.4)(44.675,29.5)(53.35,31.6)
\put(36,27.4){\circle*{0.4}}
\qbezier[60](36,27.4)(38.1,21.6)(40.2,16.8)
\put(40.2,16.8){\circle*{0.4}}
\qbezier[60](40.2,16.8)(31.65,18.6)(23.1,20.4)
\put(23.1,20.4){\circle*{0.4}}
\qbezier[58](32.6,11.2)(27.85,15.8)(23.1,20.4)
\put(32.6,11.2){\circle*{0.4}}
\qbezier[56](32.6,11.2)(25,14.3)(17.4,17.4)
\put(17.4,17.4){\circle*{0.375}}
\qbezier[55](32.6,10.8)(25,14.05)(17.4,17.3)
\put(32.6,10.8){\circle*{0.375}}
\qbezier[53](32.6,10.8)(27.85,15.25)(23.1,19.7)
\put(23.1,19.8){\circle*{0.35}}
\qbezier[50](40.2,15.6)(31.65,17.7)(23.1,19.8)
\put(40.2,15.6){\circle*{0.325}}
\qbezier[48](36,26)(38.1,20.8)(40.2,15.6)
\put(36,26){\circle*{0.3}}
\qbezier[46](36,26)(44.675,27.9)(53.35,29.8)
\put(53.35,29.8){\circle*{0.3}}
\qbezier[44](70.74,39.6)(62.045,34.7)(53.35,29.8)
\put(70.74,39.6){\circle*{0.3}}
\qbezier[40](70.74,39.6)(74.17,44.6)(77.6,49.6)
\put(77.6,49.6){\circle*{0.3}}
\qbezier[37](96,47.4)(86.4,48.5)(77.6,49.6)
\put(96,47.4){\circle*{0.3}}
\qbezier[32](96,47.4)(92.7,52.2)(89.4,57)
\put(89.4,57){\circle*{0.3}}
\qbezier[28](105,51.4)(97.2,54.5)(89.4,57)
\put(105,51.4){\circle*{0.27}}
\qbezier[23](105,51.4)(99.4,55.6)(93.8,59.8)
\put(93.8,59.8){\circle*{0.25}}

% Podstawa Sprężynki
\thinlines

\put(14.4,13){\circle*{1}}\put(2,14.6){\footnotesize{$(d_1,-\infty)$}}
\qbezier(14.4,13)(14.4,13)(27.4,5.4)
\put(27.4,5.4){\circle*{1}}\put(22,3){\footnotesize{$(0,-\infty)$}}

\qbezier[40](15.6,13.7)(22,9.75)(28.4,5.8)
\qbezier[30](15.2,13.2)(21.8,9.5)(28.4,5.8)

\put(15.7,13.7){\circle*{0.5}}
\qbezier(29.4,6.6)(29.4,6.6)(15.7,13.7)
\put(29.4,6.6){\circle*{0.6}}
\qbezier(29.4,6.6)(29.4,6.6)(17.4,14.7)
\put(17.4,14.7){\circle*{0.7}}
\qbezier(32.6,8.4)(32.6,8.4)(17.4,14.7)
\put(32.6,8.4){\circle*{0.8}}
\qbezier(32.6,8.4)(32.6,8.4)(23.1,17.4)
\put(23.1,17.4){\circle*{0.9}}
\qbezier(40.2,13.5)(40.2,13.5)(23.1,17.4)
\put(40.2,13.5){\circle*{1}}
\qbezier(40.2,13.5)(40.2,13.5)(36,24)
\put(36,24){\circle*{1}}\qbezier[48](36,24)(36,64)(36,154)
\qbezier(53.6,28.1)(53.6,28.1)(36,24)
\put(53.6,28.1){\circle*{1}}\put(51,24){\footnotesize{$(0,\infty)$}} \qbezier[120](53.6,28.1)(53.6,78.1)(53.6,158.1)
\qbezier(53.6,28.1)(53.6,28.1)(70.7,37.9)
\put(70.7,37.9){\circle*{1}}\put(68.2,33.8){\footnotesize{$(\rho,\infty)$}}\qbezier[102](70.7,37)(70.7,87)(70.7,167)
\qbezier(77.6,48)(77.6,48)(70.7,37.9)
\put(77.6,48){\circle*{1}}\put(74.6,43){\footnotesize{$(\rho_1,\infty)$}}\qbezier[82](77,48)(77,92)(77,170)
\qbezier(77.6,48)(77.6,48)(96.1,45.8)
\put(96.1,45.8){\circle*{1}}
\qbezier(89.4,55.7)(89.4,55.7)(96.1,45.8)\qbezier[42](96.1,45.8)(96.1,70.8)(96.1,160.8)
\put(89.4,55.7){\circle*{0.9}}
\qbezier(89.4,55.7)(89.4,55.7)(105,50.2)
\put(105,50.2){\circle*{0.8}}
\qbezier(93.8,58.7)(93.8,58.7)(105,50.2)
\put(93.8,58.7){\circle*{0.7}}
\qbezier(93.8,58.7)(93.8,58.7)(107.7,51.8)
\put(107.7,51.8){\circle*{0.6}}
\qbezier(95.5,59.8)(95.5,59.8)(107.7,51.8)
\put(95.5,59.8){\circle*{0.5}}

\qbezier[35](95.5,59.8)(102.05,56)(108.6,52.2)
\qbezier[25](96.1,60)(102.35,56.1)(108.6,52.2)

\put(109.5,53.1){\circle*{1}}\put(108.6,49){\footnotesize{$(\rho,$+$\infty)$}}
\qbezier(96.8,60.5)(96.8,60.5)(109.5,53.1)
\put(96.8,60.5){\circle*{1}}\put(92,63){\footnotesize{$(\rho_1,$+$\infty)$}}

% UCS
\thicklines
\put(100,4){\circle*{0.7}}
\put(100,4){\vector(3,2){10}}\put(109.2,6.8){\footnotesize{$x$}}
\put(100,4){\vector(-3,2){10}}\put(88,7){\footnotesize{$y$}}
\put(100,4){\vector(0,1){14}}\put(101.6,16){\footnotesize{$z$}}
\put(99,0){\footnotesize{$O$}}

\put(50,-8){\small{Fig. 9.}}
\put(-15,-14){\small{The inverse limit space $\X$ associated with $([0,1],f)$ in case
$f$ is $\F_{2^2}$ mapping }}
 \put(-15,-18.5){\small{possessing $\,$two intervals of$\,$ periodic points.$\,$ The projection$\,$
of $\X\,$  onto the $\,$plain}}
 \put(-15,-23){\small{$Oxy\,$ is$\,$ visualized  on$\,$ Figure$\,$ 8,$\,$ and the projection of
$\X\,$ onto the plain $\,Ozx\,$ is }}
\put(-15,-27.5){\small{visualized on Figure 7.}}
\end{picture}
\end{center}

\noindent Then the mapping $\delta$ given by
$$
\delta(a)=UaU^*, \qquad \quad a\in\A,
$$

\noindent is a unital injective endomorphism of the $C^*$-algebra $\A$, but
not an automorphism ($\delta$ is not surjective). Suppose now we
would like to investigate the  $C^*$-algebra $C^*(\A,U)$ generated
by $\A$ and $U$. If $\delta$ was an automorphism the algebra
$C^*(\A,U)$ would have the structure of the crossed-product
$\A\rtimes_{\delta} \Z$, and such objects  are very well studied.
However, we may reduce the case under consideration to the above
one by passing to a bigger $C^*$-algebra $\TA$ generated by
$\bigcup_{n\in\N} U^{*n}\A U^{n}$:
$$
\TA=C^*(\bigcup_{n\in\N} U^{*n}\A U^{n}).
$$

Indeed, cf. \cite{lebiedodzij}, then $\TA$ is a commutative
$C^*$-algebra satisfying the following relations
$$
U\TA U^*=\TA \,\,\,\, \textrm{ and }\,\,\,\, U^* \TA U= \TA.
$$

\noindent In other words mapping $\delta(\cdot)= U(\cdot)U^*$ is an
automorphism of $\TA$ and thus the algebra $C^*(\A,U)=C^*(\TA,U)$
has the structure of the crossed product $\TA\rtimes_{\delta} \Z$.
The algebra $\TA$ is called a \emph{coefficient algebra}
\cite{lebiedodzij}, \cite{maxid}, \cite{covalg} as their elements
play the role of coefficients in $C^*(\A,U)=C^*(\TA,U)$.

The question now is how the pair $(\TA,\delta)$ is related to the pair $(\A,\delta)$.
 The answer may be given in purely topological terms.

Let $X$ be the maximal ideal space  of $\A$. We recall that $X$ (also called the
spectrum of $\A$)
 is a  compact topological space, and $\A$ is isomorphic to the algebra $C(X)$ of
 continuous
 complex functions on $X$. Analogously, let $\X$ be the spectrum of $\TA$, that is
 $\TA\cong C(\X)$.

By duality it follows that  $\delta$ generates on $X$  a continuous surjective, but
not injective
mapping $f:X\to X$, and on $\X$ a homeomorphism $\f:\X\to \X$. Thus the pairs
$(\A,\delta)$ and $(\TA,\delta)$
 are determined by systems
$(X,f)$ and $(\X,\f)$ respectively, and we may reformulate our question in the
following way:

How is the system $(\X,\f)$ related to the system $(X,f)$?

The results of \cite{maxid}, \cite{covalg} give answer to this questions:
 The space $\X$ is an inverse limit of  the inverse sequence
 $X\stackrel{f}{\longleftarrow}X
\stackrel{f}{\longleftarrow} ...$, and $\f$ is the shift map induced by $f$.
In particular, inverse
limit spaces arise as spectra of coefficient algebras for crossed-products.

We note that the results of \cite{maxid} concern more general
situation than the one described above and as a consequence they
contribute to the definitions of the space $\X$
and the mapping $\f$ which generalize the inverse limit space
(\ref{inverse space}) and the induced  homeomorphism
(\ref{inverse mapping}) to the case $f$ is not necessarily
a surjective full map, see \cite{covalg}.

%\begin{center}

%\includegraphics[width=85mm]{fig3.eps}

%\vspace{5mm}

%\noindent {\small xxx}

%\end{center}

\vspace{0.5cm}

\noindent{\small Institute of Mathematics}

\noindent {\small University  in Białystok}

\noindent {\small Akademicka 2, PL-15-424 Białystok}

\noindent {\small Poland}

\noindent {\small e-mail: bartoszk@math.uwb.edu.pl}

\vspace{0.5cm}

\noindent Presented by Julian Ławrynowicz at the
Session of the Mathe\-m\-atical-Physical Commission of the \L \'od\'z
Society of Sciences and Arts on December 15, 2005

\vspace{0.5cm}

\newpage

\noindent {\bf GRANICE ODWROTNE ZWIĄZANE Z
   UNIMODALNYM\\  ODWZORWANIEM
   TYPU  $\F_{2^n}$  Z ZEROWYM SZWARCJANEM}

\vspace{0.2cm}

\noindent {\small S t r e s z c z e n i e}

{\small Praca zawiera przykład jawnego opisu granicy odwrotnej $\X$ i
jej indukowanego ho\-meo\-mor\-fizmu $\f$ przy jednym odwzorowaniu wiążącym
$f$, gdzie $f$ jest unimodalnym prze\-kształ\-ce\-niem odcinka, którego
wykres składa się z dwóch hiperbol.  Przykład ten obrazuje  zmianę
topologii $\X$ i dynamiki $\f$ wraz ze zmianą dynamiki $f$.

Ponadto praca zawiera krótki
opis  zastosowania systemów typu $(\X,\f)$ w teorii produktów krzyżowych
$C^*$-algebr. }

%\begin{thebibliography}{99}
%{caban}P.~Caban, bleble

%\end{thebibliography}

\end{document}